\documentclass[11pt,letterpaper]{amsart}
\pdfoutput=1

\usepackage[normalem]{ulem}

\makeatletter
\newcommand*{\inlineequation}[2][]{%
  \begingroup
    \refstepcounter{equation}%
    \ifx\\#1\\%
    \else
      \label{#1}%
    \fi
    \relpenalty=10000 %
    \binoppenalty=10000 %
    \ensuremath{%
      #2%
    }%
    ~\@eqnnum
  \endgroup
}
\makeatother

\usepackage{tikz-cd}

\usepackage{tcolorbox}

\usepackage[arrow,curve,matrix,tips,2cell]{xy}

\usepackage[hmargin=2.5cm,vmargin=1.5cm]{geometry}

\usepackage{bm}
\usepackage{fancyhdr}

\usepackage{amsmath}
\usepackage{amscd}
\usepackage{amssymb}
\usepackage{latexsym}

\usepackage{graphicx}

\usepackage{mathrsfs}

\usepackage[shortlabels]{enumitem}
\usepackage{hyperref}
\usepackage{comment}

\newtheorem{theorem}{Theorem}[section]
\newtheorem*{theorem*}{Theorem}
\newtheorem{lemma}[theorem]{Lemma}
\newtheorem*{lemma*}{Lemma}
\newtheorem{corollary}[theorem]{Corollary}
\newtheorem*{corollary*}{Corollary}
\newtheorem{proposition}[theorem]{Proposition}

\newtheorem{remark}[theorem]{Remark}
\newtheorem{definition}[theorem]{Definition}
\newtheorem*{definition*}{Definition}

\newtheorem{question}[theorem]{Question}
\newtheorem*{question*}{Question}

\newtheorem{example}[theorem]{Example}
\newtheorem{examples}[theorem]{Examples}

\newtheorem{thm}{Theorem}[subsection]

\newtheorem{cor}[thm]{Corollary}

\makeatletter
\def\revddots{\mathinner{\mkern1mu\raise\p@
\vbox{\kern7\p@\hbox{.}}\mkern2mu
\raise4\p@\hbox{.}\mkern2mu\raise7\p@\hbox{.}\mkern1mu}}
\makeatother 
\newcommand{\bgl}{\begin{equation}} 
\newcommand{\egl}{\end{equation}}
\newcommand{\bgloz}{\begin{equation*}} 
\newcommand{\egloz}{\end{equation*}}
\newcommand{\bgln}{\begin{eqnarray}} 
\newcommand{\egln}{\end{eqnarray}}
\newcommand{\bglnoz}{\begin{eqnarray*}} 
\newcommand{\eglnoz}{\end{eqnarray*}}
\newcommand{\btheo}{\begin{theorem}}
\newcommand{\etheo}{\end{theorem}}
\newcommand{\btheooz}{\begin{theorem*}}
\newcommand{\etheooz}{\end{theorem*}}

\newcommand{\blemma}{\begin{lemma}}
\newcommand{\elemma}{\end{lemma}}
\newcommand{\blemmaoz}{\begin{lemma*}}
\newcommand{\elemmaoz}{\end{lemma*}}
\newcommand{\bproof}{\begin{proof}}
\newcommand{\eproof}{\end{proof}}
\newcommand{\bbew}{\begin{beweis}}
\newcommand{\ebew}{\end{beweis}}
\newcommand{\bremark}{\begin{remark}}
\newcommand{\eremark}{\end{remark}}
\newcommand{\bdefin}{\begin{definition}}
\newcommand{\edefin}{\end{definition}}
\newcommand{\bdefinoz}{\begin{definition*}}
\newcommand{\edefinoz}{\end{definition*}}
\newcommand{\bex}{\begin{example}}
\newcommand{\eex}{\end{example}}
\newcommand{\bexs}{\begin{examples}}
\newcommand{\eexs}{\end{examples}}

\newcommand{\bprop}{\begin{proposition}}
\newcommand{\eprop}{\end{proposition}}
\newcommand{\bcor}{\begin{corollary}}
\newcommand{\ecor}{\end{corollary}}
\newcommand{\bfa}{\begin{cases}} 
\newcommand{\efa}{\end{cases}}

\newcommand{\bquestion}{\begin{question}}
\newcommand{\equestion}{\end{question}}
\newcommand{\bquestionoz}{\begin{question*}}
\newcommand{\equestionoz}{\end{question*}}
\newcommand{\cA}{\mathcal A}

\newcommand{\cC}{\mathcal C}

\newcommand{\cE}{\mathcal E}

\newcommand{\cG}{\mathcal G}

\newcommand{\cL}{\mathcal L}
\newcommand{\cM}{\mathcal M}

\newcommand{\cO}{\mathcal O}
\newcommand{\cP}{\mathcal P}
\newcommand{\cQ}{\mathcal Q}

\newcommand{\cS}{\mathcal S}
\newcommand{\cT}{\mathcal T}

\newcommand{\cW}{\mathcal W}

\def\Cz{\mathbb{C}}

\def\Iz{\mathbb{I}}

\def\Lz{\mathbb{L}}

\def\Nz{\mathbb{N}}

\def\Qz{\mathbb{Q}}

\def\Zz{\mathbb{Z}}

\def\1z{\mathbb{1}}
\newcommand{\fA}{\mathfrak A}
\newcommand{\fB}{\mathfrak B}
\newcommand{\fC}{\mathfrak C}
\newcommand{\fD}{\mathfrak D}

\newcommand{\fN}{\mathfrak N}

\newcommand{\fT}{\mathfrak T}

\newcommand{\ti}{\tilde}

\newcommand{\ma}{\mapsto} 
\newcommand{\onto}{\twoheadrightarrow} 
\newcommand{\into}{\hookrightarrow} 

\def\SEMI{\mbox{$\times\kern-2pt\vrule height5pt width.6pt \kern3pt $}}

\newcommand{\Hom}{{\rm Hom}\,}
\newcommand{\End}{{\rm End}\,}
\newcommand{\Aut}{{\rm Aut}}

\newcommand{\Sp}{{\rm Sp\,}}

\newcommand{\id}{{\rm id}}

\renewcommand{\dim}{{\rm dim}}
\newcommand{\rk}{{\rm rk}}

\renewcommand{\ker}{{\rm ker}\,}

\newcommand{\lcm}{{\rm lcm}} 
\newcommand{\reg}{^\times} 
\newcommand{\lspan}{{\rm span}} 
\newcommand{\defeq}{\mathrel{:=}} 

\newcommand{\dop}{\text{: }} 
\newcommand{\ilim}{\varinjlim} 
\newcommand{\plim}{\varprojlim} 
\newcommand{\supp}{{\rm supp}} 
\newcommand{\lge}{\left\{} 
\newcommand{\rge}{\right\}} 
\newcommand{\lsp}{\left\langle} 
\newcommand{\rsp}{\right\rangle} 
\newcommand{\gekl}[1]{\lge #1 \rge} 
\newcommand{\spkl}[1]{\lsp #1 \rsp} 
\newcommand{\menge}[2]{\gekl{ #1 \dop #2 }} 
\newcommand{\mfa}{\mathfrak{a}}
\newcommand{\mfb}{\mathfrak{b}}
\newcommand{\mfc}{\mathfrak{c}}
\newcommand{\mfg}{\mathfrak{g}}
\newcommand{\mfk}{\mathfrak{k}}
\newcommand{\mfl}{\mathfrak{l}}

\newcommand{\mfs}{\mathfrak{s}}
\newcommand{\mft}{\mathfrak{t}}
\newcommand{\scA}{\mathscr{A}}
\newcommand{\scB}{\mathscr{B}}
\newcommand{\scS}{\mathscr{S}}
\renewcommand{\k}{\mathfrak{k}}

\newcommand{\p}{\mathfrak{p}}

\newcommand{\m}{\mathfrak{m}}

\newcommand{\Gal}{\textup{Gal}}

\newcommand{\acts}{\curvearrowright}

\def\tor{{\rm tor}}

\newcommand{\ord}{\textup{ord}} 
\newcommand{\spn}{\textup{span}}

\newcommand{\ol}[1]{\overline{#1}}
\newcommand{\GL}{\textup{GL}}
\newcommand{\SL}{\textup{SL}}
\newcommand{\M}{\textup{M}}
\newcommand{\Mult}{\textup{Mult}}

\newcommand{\gp}[1]{\langle #1\rangle}
\newcommand{\mon}[1]{\langle #1\rangle^+}

\newcommand{\im}{\textup{im}}
\newcommand{\msc}{\mathscr}
\newcommand{\bd}{\partial\widehat{\cE}}
\renewcommand{\sp}{\textup{sp}}

\begin{document}

\title{Algebraic actions II. Groupoid rigidity}

\thispagestyle{fancy}

\author{Chris Bruce}
\thanks{C. Bruce was supported by a Banting Fellowship administered by the Natural Sciences and Engineering Research Council of Canada (NSERC) and has received funding from the European Union’s Horizon 2020 research and innovation programme under the Marie Sklodowska-Curie grant agreement No 101022531. X. Li has received funding from the European Research Council (ERC) under the European Union’s Horizon 2020 research and innovation programme (grant agreement No. 817597).}
\address{School of Mathematics, Statistics and Physics, Herschel Building, Newcastle University, Newcastle upon Tyne, NE1 7RU, United Kingdom}
\email{chris.bruce@newcastle.ac.uk}

\author{Xin Li}
\address{School of Mathematics and Statistics, University of Glasgow, University Place, Glasgow G12 8QQ, United Kingdom}
\email{Xin.Li@glasgow.ac.uk}



\begin{abstract}
We establish several surprising rigidity results for \'etale groupoids associated with algebraic actions.
Our initial motivation comes from algebraic number theory, where our results are witnessed in a particularly striking fashion: We prove that the groupoids associated with the action of the multiplicative monoid in a ring of algebraic integers on the additive group of the ring remembers the ring up to isomorphism. Already this special case resolves an open problem about isomorphisms of Cartan pairs and leads to a dynamical analogue of the Neukirch--Uchida theorem using topological full groups.
\end{abstract}


\maketitle

\section{Introduction}
\setlength{\parindent}{0cm} \setlength{\parskip}{0cm}

The goal of this paper is to establish rigidity results for groupoids arising from algebraic actions (see \cite{BL2}). In the special situation of actions coming from algebraic number theory, which was our original motivation, our main result takes the following form: Let $K_1$ and $K_2$ be algebraic number fields with rings of algebraic integers $R_1$ and $R_2$, respectively. The $ax+b$-group $K_i\rtimes K_i^\times$ acts on the compact ring $\ol{R_i}$ of integral adeles over $K_i$, so we may form the partial transformation groupoid $(K_i\rtimes K_i^\times)\ltimes \ol{R_i}$. 
\begin{theorem*}[Corollary~\ref{cor:IntroNT}]
The fields $K_1$ and $K_2$ are isomorphic if and only if $(K_1\rtimes K_1^\times)\ltimes \ol{R_1}$ and $(K_2\rtimes K_2^\times)\ltimes \ol{R_2}$ are isomorphic as topological groupoids.
\end{theorem*}

Let us explain two applications of this theorem. The first is to Cartan pairs in semigroup C*-algebras. The C*-algebra of $(K_i\rtimes K_i^\times)\ltimes \ol{R_i}$ is canonically isomorphic to the C*-algebra $C^*(R_i\rtimes R_i^\times)$ of the $ax+b$-semigroup $R_i\rtimes R_i^\times$ as considered in \cite{CDL,CEL,Li:Adv2014,Li:Adv2016}, where the multiplicative monoid $R_i^\times:=R_i\setminus\{0\}$ acts on (the additive group of) $R_i$ by multiplication. Using this, we obtain the following result, which resolves a natural problem left open in \cite{Li:Adv2016} (we even resolve the more general problem left open in \cite{BL}, see Remark~\ref{rmk:resolve}).
\begin{corollary*}[Corollary~\ref{cor:IntroNTCartan}]
The fields $K_1$ and $K_2$ are isomorphic if and only if there is a Cartan-preserving *-isomorphism between $C^*(R_1\rtimes R_1^\times)$ and $C^*(R_2\rtimes R_2^\times)$.
\end{corollary*}
\setlength{\parindent}{0cm} \setlength{\parskip}{0cm}

The second application we wish to mention here is a dynamical analogue of the Neukirch--Uchida theorem using topological full groups. Each groupoid $(K_i\rtimes K_i^\times)\ltimes \ol{R_i}$ gives rise to a topological full group $\bm{F}((K_i\rtimes K_i\reg)\ltimes\ol{R_i})$ as defined in \cite{Mat12,Nek19}, which is a countable discrete group of homeomorphisms of the Cantor space $\ol{R_i}$.
\begin{corollary*}[Corollary~\ref{cor:IntroNT}]
The fields $K_1$ and $K_2$ are isomorphic if and only if the topological full groups $\bm{F}((K_1\rtimes K_1\reg)\ltimes\ol{R_1})$ and $\bm{F}((K_2\rtimes K_2\reg)\ltimes\ol{R_2})$ are isomorphic as abstract groups.
\end{corollary*}
\setlength{\parindent}{0cm} \setlength{\parskip}{0cm}

The present paper contains various generalizations and versions of the above rigidity results for partial transformation groupoids associated with general algebraic actions of semigroups in \cite{BL2}: If two such groupoids (satisfying appropriate conditions) are isomorphic, then the globalizations of the initial algebraic actions rationally embed in each other. For specific example classes arising for instance from toral endomorphisms, actions from rings, or actions from commutative algebra, this mutual embedability can be improved in various ways to obtain surprisingly strong rigidity phenomena, as illustrated above in the number-theoretic setting.

\subsection*{Context}
Over the last few decades, there have been significant advances on interactions between groups, C*-algebras, and topological dynamics. 
C*-algebras of dynamical origin, such as crossed product and groupoid C*-algebras, play a central role in operator algebras, and there are many deep connections between operator algebras and topological dynamics, for instance, the study of K-theoretic invariants of C*-algebras associated with Cantor minimal systems led to striking classification results for orbit equivalence by Giordano, Matui, Putnam, and Skau \cite{GPS1,GMPS1,GMPS2}. At the heart of the interplay between topological dynamics and C*-algebras are the notions of \emph{\'{e}tale groupoids} (cf. \cite{Ren80} and \cite[\S~II]{SSW20}) and \emph{Cartan subalgebras} (see \cite{Kum:86,Ren}). The bridge to group theory is built by \emph{Topological full groups} attached to \'{e}tale groupoids, which were considered by Giordano, Putnam, and Skau in the setting of Cantor minimal systems \cite{GPS2} and in general by Matui \cite{Mat15}. They are defined as groups of homeomorphisms of the Cantor set whose `graph' is contained in the groupoid as a compact open subset. These groups are very interesting in their own right, and they have given rise to several celebrated example classes, resolving open problems in group theory  \cite{JM,JNdlS,Nek18,SWZ}.

In many cases, the groupoid, the topological full group, and the Cartan pair encode exactly the same information, and there is a precise dictionary how to translate between these---a priori very different---languages. This allows one to bring tools from one area to bear on problems in another, e.g., groupoid techniques are a key tool in understanding topological full groups, and the language of Cartan pairs was crucial for the result in \cite{LiInvent} that all classifiable C*-algebras have groupoid models.

\subsection*{Algebraic actions} 
A particularly interesting arena where aforementioned interactions play out is in the context of C*-algebras and groupoids associated with \emph{algebraic actions of semigroup}. Such actions are the irreversible counterparts of \emph{algebraic actions of groups} (cf. \cite[Chapters~13~and~14]{KerrLi} and \cite{Sch}), which have been studied intensely since the late 80s. In contrast to the reversible case, the theory of algebraic actions of semigroups is much less developed.

In \cite{BL2}, we initiated the study of algebraic actions of semigroups from the point of view of \'etale groupoids. On the one hand, the example class of algebraic actions is interesting and exhibits new phenomena, and on the other hand, due to the particular structure of algebraic actions, a variety of tools is available, allowing for a systematic and detailed study. An interesting new phenomenon that arises in this new setting is that actions by non-invertible endomorphisms of a given group automatically produce a particular completion of the group, and the original action induces a system of partial homeomorphisms on this completion. Using this, we associated an \'e{tale} groupoid to each algebraic action of a semigroup \cite{BL2}; these groupoids are interesting in their own right but also give access to analyzing properties of C*-algebras generated by natural representations of the initial algebraic action, which was our goal in \cite{BL2}. Our goal here---with applications to topological full groups in mind---is to study the natural question of how much information the groupoids constructed in \cite{BL2} remember about the original algebraic actions. Surprisingly, we discover the phenomenon of \emph{groupoid rigidity} for a variety of example classes, i.e., our groupoids remember more information than expected -- in special cases, they even remember everything. In the number-theoretic setting as mentioned above, our rigidity results resolve an open problem about Cartan subalgebras for the C*-algebras of $ax+b$-semigroups introduced by Cuntz, Deninger, and Laca in \cite{CDL} and lead to a dynamical analogue of the celebrated Neukirch--Uchida theorem. However, we wish to stress that these are very special cases of our general theorems, and only serve to demonstrate the surprising groupoid rigidity phenomenon arising in the setting of general algebraic actions.

In our setting, the groupoids are completely characterized by their topological full groups, so it becomes particularly interesting to study properties of these groups, and advances in this setting has the potential for a high impact both in topological dynamics and in group theory. This is for example the subject of forthcoming work of the first-named author with Y. Kubota and T. Takeishi.

\subsection*{Results}
Let us now formulate our rigidity results. An algebraic action $\sigma \colon S \acts A$ consists of a monoid $S$, an Abelian group $A$, and a semigroup homomorphism from $S$ to injective endomorphisms of $A$, denoted by $S \to \End(A), \, s \ma \sigma_s$. We will always assume our algebraic actions to be non-automorphic (i.e., not all $\sigma_s$ are automorphisms) and faithful (i.e., the map $s \ma \sigma_s$ is injective). Let $\cC$ be the collection of subgroups of $A$ which are of the form $\sigma_{t_1}^{-1} \sigma_{s_1} \dotsm \sigma_{t_m}^{-1} \sigma_{s_m} A$, where $\sigma_t^{-1}X \defeq \menge{a \in A}{\sigma_t(a) \in X}$ for a subset $X \subseteq A$. In this paper, we will always assume that $\sigma$ has the finite index property \eqref{FI}, i.e., $\# (A/C) < \infty$ for all $C \in \cC$, or equivalently, $\# (A / \sigma_sA) < \infty$ for all $s \in S$. In this case, the completion of $A$ mentioned above is given by $\ol{A} \defeq \plim_{C \in \cC} A/C$, where $\cC$ is partially ordered by inclusion. The standing assumptions in this paper will furthermore include that $\sigma \colon S \acts A$ admits a globalization (which is always assumed to be minimal in the sense of Remark~\ref{rem:mscA=CUPA}), i.e., an embedding of $S$ into a group $\msc{S}$ and a group $\msc{A}$ containing $A$ together with an algebraic action $\tilde{\sigma} \colon \msc{S} \acts \msc{A}$ (necessarily by automorphisms) such that $\tilde{\sigma}_s \vert_A = \sigma_s$ for all $s\in S$. This allows us to form a partial action of $\msc{S}$ on $A$ by letting $s \in \msc{S}$ act via the restriction $A \cap \ti{\sigma}_s^{-1}A \to \ti{\sigma}_s A \cap A$ of $\ti{\sigma}_s$ to $A \cap \ti{\sigma}_s^{-1}A$. Similarly, we also have the partial action of $\msc{A}$ on $A$ by translation. In \cite{BL2}, we identified a condition, called \eqref{JF}, which ensures that these partial actions on $A$ extend to a partial action of $\msc{A} \rtimes \msc{S}$ on $\ol{A}$ by partial homeomorphisms (see \S~\ref{ss:Standing} and \cite[\S~3.3]{BL2} for details). In this paper, we will always assume that \eqref{JF} holds. The associated partial transformation groupoid $\cG_{\sigma} \defeq (\msc{A} \rtimes \msc{S}) \ltimes \ol{A}$ is the groupoid constructed in \cite[\S~3]{BL2} arising from our algebraic action $\sigma$. We can now state our main rigidity result.
\begin{thm}[see Theorem~\ref{thm:frmain}]
Let $\sigma \colon S \acts A$ and $\tau \colon T \acts B$ be two algebraic actions of monoids as above, with globalizations $\tilde{\sigma} \colon \msc{S} \acts \msc{A}$, $\tilde{\tau} \colon \msc{T} \acts \msc{B}$ and groupoids $\cG_{\sigma}$, $\cG_{\tau}$. Assume that $S$ and $T$ are Abelian, that $\msc{A}$  and $\msc{B}$ are torsion-free and finite rank, and that there exist $s \in S$ and $t \in T$ such that $\id_A - \sigma_s \colon A \to A$ and $\id_B - \tau_t \colon B \to B$ are injective. 

If $\cG_{\sigma}$ and $\cG_{\tau}$ are isomorphic as topological groupoids, then $\tilde{\sigma} \colon \msc{S} \acts \msc{A}$ and $\tilde{\tau} \colon \msc{T} \acts \msc{B}$ embed rationally into each other, i.e., there exist injective homomorphisms $\mft \colon \msc{S} \into \msc{T}$, $\mfb \colon \Qz \otimes \msc{A} \into \Qz \otimes \msc{B}$, $\mfs \colon \msc{T} \into \msc{S}$, and $\mfa \colon \Qz \otimes \msc{B} \into \Qz \otimes \msc{A}$ such that $\mfb((\id_{\Qz} \otimes \ti{\sigma}_s)(x)) = (\id_{\Qz} \otimes \ti{\tau}_{\mft(s)})(\mfb(x))$ and $\mfa((\id_{\Qz} \otimes \ti{\tau}_t)(y)) = (\id_{\Qz} \otimes \ti{\sigma}_{\mfs(t)})(\mfa(y))$ for all $s \in \msc{S}$, $x \in \Qz \otimes \msc{A}$, $t \in \msc{T}$, and $y \in \Qz \otimes \msc{B}$.
\end{thm}
\setlength{\parindent}{0cm} \setlength{\parskip}{0cm}

We refer the reader to \S~\ref{s:c}, in particular \S~\ref{ss:Cons}, for more general results and further explanations.
\setlength{\parindent}{0cm} \setlength{\parskip}{0.5cm}

Let us now present a first class of algebraic actions where our general rigidity result applies.
\begin{cor}[see Corollary~\ref{cor:mixing}]
Assume that $S$ and $T$ are Abelian, torsion-free monoids, that $A$ and $B$ are torsion-free Abelian groups of finite rank, and that $\sigma \colon S \acts A$ and $\tau \colon T \acts B$ are non-automorphic faithful algebraic actions. 
Further suppose that there exist $s\in S$ and $t\in T$ such that the dual actions $\hat{\sigma}\vert_{s^\Nz}:s^\Nz\acts\widehat{A}$ and $\hat{\tau}\vert_{t^\Nz}:t^\Nz\acts\widehat{B}$ are mixing, where $s^\Nz:=\{s^n : n\in\Nz\}$ and $t^\Nz:=\{t^n : n\in\Nz\}$.
Let $\ti{\sigma} \colon S^{-1} S \acts S^{-1} A$ and $\ti{\tau} \colon T^{-1} T \acts T^{-1} B$ be the canonical globalizations as in \cite[Example~2.4]{BL2}, and denote the groupoids attached to $\sigma$ and $\tau$ by $\cG_{\sigma}$ and $\cG_{\tau}$, respectively.

If $\cG_{\sigma}$ and $\cG_{\tau}$ are isomorphic as topological groupoids, then there exist injective homomorphisms $\mft\colon S^{-1} S \hookrightarrow T^{-1} T$ and $\mfb \colon S^{-1} A \hookrightarrow  T^{-1} B$ such that $\mfb(\ti{\sigma}_s(x))=\ti{\tau}_{\mft(s)}(\mfb(x))$ for all $s \in S^{-1} S$ and $x\in S^{-1} A$, injective homomorphisms $\mfs\colon T^{-1} T \hookrightarrow S^{-1} S$ and $\mfa \colon T^{-1} B \hookrightarrow  S^{-1} A$ such that $\mfa(\ti{\tau}_t(y))=\ti{\sigma}_{\mfs(t)}(\mfa(y))$ for all $t \in T^{-1} T$ and $y \in T^{-1} B$, and the images of $\mfb$ and $\mfa$ are finite index subgroups of $T^{-1} B$ and $S^{-1}A$, respectively.
\end{cor}
\setlength{\parindent}{0cm} \setlength{\parskip}{0cm}

The reader may consult \S~\ref{ss:mixing} for more explanations and details. The concrete case of toral endomorphisms is treated in Example~\ref{ex:ToralEndo}.
\setlength{\parindent}{0cm} \setlength{\parskip}{0.5cm}

Another motivating example class is given by the action of the monoid of non-zerodivisors of a ring on the additive group of the ring by multiplication. For instance, torsion-free commutative rings which are finitely generated as additive groups have received a great deal of attention because of Bhargava's work \cite{Bh1,Bh2,Bh3,Bhd1,Bh4,Bhd2}. In this setting, our rigidity result implies the following.
\begin{thm}[see Theorem~\ref{thm:Bhargava}]
\label{thm:IntroCommRing}
Let $R_i$, $i=1,2$, be finitely generated torsion-free commutative rings. For $i=1,2$, let $R_i\reg$ be the monoid of non-zerodivisors in $R_i$ and $\sigma_i \colon R_i\reg \acts R_i$ the algebraic action given by multiplication. If the corresponding groupoids $\cG_{\sigma_1}$ and $\cG_{\sigma_2}$ are isomorphic, then $\Qz \otimes R_1$ and $\Qz \otimes R_2$ are isomorphic as $\Qz$-algebras.
\end{thm}

Somewhat surprisingly, using very different methods, we obtain a similar rigidity result for special classes of non-commutative rings. 
\begin{thm}[see Corollary~\ref{cor:ssQ}]
\label{thm:IntroNCRing}
For $i=1,2$, let $R_i$ be a ring whose additive group is finitely generated and torsion-free. Let $R_i\reg$ be the monoid of left regular elements (i.e., non-left-zerodivisors) in $R_i$ and $\sigma_i \colon R_i\reg \acts R_i$ the algebraic action given by left multiplication. Suppose that $\Qz\otimes R_1$ and $\Qz\otimes R_2$ are semisimple $\Qz$-algebras. If the corresponding groupoids $\cG_{\sigma_1}$ and $\cG_{\sigma_2}$ are isomorphic, then $\Qz\otimes R_1$ and $\Qz\otimes R_2$ are isomorphic as $\Qz$-algebras.
\end{thm}
\setlength{\parindent}{0cm} \setlength{\parskip}{0cm}

Examples of rings that are covered by Theorem~\ref{thm:IntroNCRing} include integral group rings of finite groups and rings of matrices over orders in algebraic number fields.
\setlength{\parindent}{0cm} \setlength{\parskip}{0.5cm}

Let us now apply Theorem~\ref{thm:IntroCommRing} to rings of algebraic integers in number fields. Let $R$ be such a ring, with quotient field $K$. Consider the algebraic action $\sigma \colon R\reg \acts R$ by multiplication. The completion $\ol{R}$ is given in this case by the (additive group of the) integral adele ring, and the partial action from above is given by the canonical partial action $K \rtimes K\reg \acts \ol{R}$. 
Moreover, the groupoid $\cG_{\sigma}$ in this case is the partial transformation groupoid $(K\rtimes K\reg)\ltimes\ol{R}$, and its C*-algebra coincides with the ring C*-algebra $\fA[R]$, which has been introduced and studied in \cite{Cuntz,CuLi,Li:Ring}.
Since $(K\rtimes K\reg)\ltimes\ol{R}$ is effective, $\fA[R]$ contains a canonical Cartan subalgebra $\fD[R]$. Moreover, consider the topological full group $\bm{F}((K\rtimes K\reg)\ltimes\ol{R})$ of $(K\rtimes K\reg)\ltimes\ol{R}$ given by the group of global bisections (see for instance \cite{Mat12,Nek19}) and its commutator subgroup $\bm{D}((K\rtimes K\reg)\ltimes\ol{R})$. 
\begin{cor}[see Corollary~\ref{cor:NT}]
\label{cor:IntroNT}
Let $R_i$, $i=1,2$, be rings of algebraic integers in number fields $K_i$. With the notation introduced above, the following are equivalent:
\setlength{\parindent}{0cm} \setlength{\parskip}{0cm}

\begin{enumerate}[\upshape(i)]
    \item $K_1$ and $K_2$ are isomorphic;
    \item $(K_1\rtimes K_1\reg)\ltimes\ol{R_1}$ and $(K_2\rtimes K_2\reg)\ltimes\ol{R_2}$ are isomorphic as topological groupoids;
    \item $K_1 \rtimes K_1\reg \acts \ol{R_1}$ and $K_2 \rtimes K_2\reg \acts \ol{R_2}$ are continuously orbit equivalent in the sense of \cite{Li:ETDS,Li:IMRN2017};
    \item $(\fA[R_1],\fD[R_1])$ and $(\fA[R_2],\fD[R_2])$ are isomorphic as Cartan pairs;
    \item $\bm{F}((K_1\rtimes K_1\reg)\ltimes\ol{R_1})$ and $\bm{F}((K_2\rtimes K_2\reg)\ltimes\ol{R_2})$ are isomorphic as abstract groups;
    \item $\bm{D}((K_1\rtimes K_1\reg)\ltimes\ol{R_1})$ and $\bm{D}((K_2\rtimes K_2\reg)\ltimes\ol{R_2})$ are isomorphic as abstract groups.
\end{enumerate}
\end{cor}
\setlength{\parindent}{0cm} \setlength{\parskip}{0.5cm}

\bremark
The equivalences of (i), (v), and (vi) in Corollary~\ref{cor:IntroNT} gives dynamical analogues of the Neukirch--Uchida theorem from anabelian geometry which says that the absolute Galois group of a number field remembers the field up to isomorphism \cite{Neu:69,Uchida}. A different dynamical version of the Neukirch--Uchida theorem is given in \cite{BT} using completely different techniques. Moreover, the structure of our groups is much different from the absolute Galois groups or the topological full groups from \cite{BT} since, e.g., $\bm{D}((K\rtimes K\reg)\ltimes\ol{R})$ is simple.
\eremark

\bremark
The equivalence of (i) and (iv) in Corollary~\ref{cor:IntroNT} is in stark contrast with \cite[Corollary~1.3]{LiLu}, which says that the ring C*-algebras $\fA[R_1]$ and $\fA[R_2]$ are always isomorphic.
One consequence of this is that $\fA[\Zz]$ contains a family, parameterized by all number fields, of isomorphic but non-conjugate Cartan subalgebras.
\eremark

Semigroup C*-algebras $C^*(R \rtimes R\reg)$ of $ax+b$-semigroups $R \rtimes R\reg$ were studied in \cite{CDL,CEL,Li:Adv2014,Li:Adv2016} for rings of algebraic integers $R$ in number fields $K$. $C^*(R \rtimes R\reg)$ has a canonical groupoid model (see \cite{CuLi,Li:Adv2014}), and hence contains a canonical Cartan subalgebra $D(R \rtimes R\reg)$. It was shown in \cite{Li:Adv2016} how to recover the Dedekind zeta function and the ideal class group of $K$ from the Cartan pair $(C^*(R \rtimes R\reg), D(R \rtimes R\reg))$. However, the natural question of whether $(C^*(R \rtimes R\reg), D(R \rtimes R\reg))$ completely determines the number field $K$ has been left open in \cite{Li:Adv2016} (see the question at the end of \cite[\S~1]{Li:Adv2016}). Since $(\fA[R],\fD[R])$ can be recovered from $(C^*(R \rtimes R\reg), D(R \rtimes R\reg))$, we are now able to answer this question.
\begin{cor}
\label{cor:IntroNTCartan}
Let $R_i$, $i=1,2$, be rings of algebraic integers in number fields $K_i$. Then $K_1 \cong K_2$ if and only if $(C^*(R_1 \rtimes R_1\reg), D(R_1 \rtimes R_1\reg))$ and $(C^*(R_2 \rtimes R_2\reg), D(R_2 \rtimes R_2\reg))$ are isomorphic as Cartan pairs.
\end{cor}
\setlength{\parindent}{0cm} \setlength{\parskip}{0cm}

In fact, we completely resolve the more general problem left open in \cite[\S~5.2]{BL}, see Remark~\ref{rmk:resolve}.
\setlength{\parindent}{0cm} \setlength{\parskip}{0.5cm}

Theorem~\ref{thm:IntroNCRing}, applied to matrix algebras over rings of algebraic integers, yields the following analogous result. 
\begin{cor}[see Example~\ref{ex:MnK}]
For $i=1,2$, let $R_i$ be the ring of integers in a number field $K_i$, and let $n_i$ be a positive integer. Consider the algebraic action $\sigma_i \colon \M_{n_i}(R_i)\reg \acts \M_{n_i}(R_i)$ by multiplication. The groupoids $\cG_{\sigma_1}$ and $\cG_{\sigma_2}$ are isomorphic if and only if $n_1=n_2$ and $K_1 \cong K_2$. 
\end{cor}
\setlength{\parindent}{0cm} \setlength{\parskip}{0cm}

We also have further equivalent statements analogous to (iii) -- (vi) in Corollary~\ref{cor:IntroNT}, and the analogue of Corollary~\ref{cor:IntroNTCartan} holds as well.
\setlength{\parindent}{0cm} \setlength{\parskip}{0.5cm}

We would like to point out that we obtain more general results than the ones presented in this introduction (see \S~\ref{s:c} for details). Thus, we can additionally treat the following example classes:
\setlength{\parindent}{0cm} \setlength{\parskip}{0cm}
\begin{enumerate}[\upshape(a)]
    \item Semigroups of canonical endomorphisms of finite rank torsion-free Abelian groups (\S~\ref{ss:endoabgroups});
    \item Actions form adding scalars to algebraic actions of subgroups of special linear groups (\S~\ref{ss:addingscalars});
    \item Arithmetical $\cS$-integer dynamical systems (\S~\ref{ss:arithmeticdynamics});
    \item Actions of congruence monoids on rings of algebraic integers (\S~\ref{ssecNT});
    \item $\Nz^d$-actions coming from zero-dimensional ideals in commutative algebra (\S~\ref{sss:commalg});
    \item Actions from integral group rings of finite groups (Example~\ref{ex:grprings});
    \item Actions from orders in central simple algebras over number fields (Example~\ref{ex:ordersinssalgs}).
\end{enumerate}
\setlength{\parindent}{0cm} \setlength{\parskip}{0.5cm}

The proofs of our rigidity results are inspired by continuous orbit equivalence rigidity for odometers (see \cite{CM}). Given an algebraic action $\sigma \colon S \acts A$ satisfying our standing assumptions, the restriction of the partial action of $\msc{A} \rtimes \msc{S} \acts \ol{A}$ to $A$ yields an odometer action $A \acts \ol{A}$. Our key insight is that a careful analysis allows us to identify situations where rigidity can be upgraded from these odometer actions to groupoid rigidity in our sense. For Abelian acting monoids, we take advantage of algebraic identities in semidirect product groups arising from our algebraic actions. For non-Abelian acting monoids, our rigidity results rely on the structure of nilpotent elements in certain matrix algebras.

\section{Preliminaries}
\label{sec:prelim}

\subsection{Standing assumptions}
\label{ss:Standing}

We first explain the standing assumptions on algebraic actions that we will assume in this paper. Let $S$ be a non-trivial left cancellative monoid and $A$ an Abelian group, written additively, with identity element $0$. Assume $\sigma\colon S\acts A$ is an \emph{algebraic $S$-action}, i.e., $\sigma\colon S\to \End_\Zz(A), \, s \ma \sigma_s$ is a monoid homomorphism such that $\sigma_s$ is an injective endomorphism $A \to A$ for all $s\in S$. Actions of this form are called \emph{algebraic actions} (cf. \cite{BL2}). Unless $S$ is a group, we shall assume that $\sigma\colon S\acts A$ is non-automorphic, i.e., there exists $s\in S$ such that $\sigma_sA\subsetneq A$. This in particular implies that $A$ is non-trivial. We will also always assume that the action $\sigma\colon S\acts A$ is \emph{faithful}, i.e., $s\mapsto\sigma_s$ is injective. 

Let $\cC=\cC_{S\acts A}$ be the family of \emph{$S$-constructible subgroups of $A$}, i.e., $\cC$ is the smallest collection of subgroups of $A$ such that 
\setlength{\parindent}{0cm} \setlength{\parskip}{0cm}
\begin{itemize}
    \item $A\in\cC$;
    \item $s\in S$ and $C\in\cC$ implies $\sigma_s C,\sigma_s^{-1}C\in\cC$.
\end{itemize}
It follows that $\cC$ is closed under taking finite intersections (see \cite[Proposition~3.9]{BL2}. We say that $\sigma\colon S\acts A$ satisfies the \emph{finite index property} (see \cite[Definition~7.1]{BL2}) if
\begin{equation}
\label{FI}\tag{FI}
	\# (A/\sigma_sA) <\infty\quad \text{for all } s\in S,
\end{equation}
If $\sigma\colon S\acts A$ satisfies \eqref{FI}, then \cite[Proposition~7.2]{BL2} implies that every member of $\cC$ is a finite index subgroup of $A$. In this case, we get a compact group $\ol{A}:=\varprojlim_{C\in\cC}A/C$, and the canonical homomorphism $A\to\ol{A}$ has kernel $\bigcap_{C\in\cC}C$. Given $C \in \cC$, we denote by $\ol{C}$ the kernel of the canonical projection $\ol{A} \onto A/C$, i.e., the preimage of $C \in A/C$ in $\ol{A}$. We have a canonical homeomorphism $\ol{C} \cong \plim_{D \in \cC, \, D \subseteq C} C/D$.
\setlength{\parindent}{0cm} \setlength{\parskip}{0.5cm}

\bremark
It suffices to check \eqref{FI} for generators of $S$: Say $S$ is generated by $\cS$. Then we can proceed inductively on the word length with respect to $\cS$ of an element in $S$. Suppose $\# (A / \sigma_s A) < \infty$ for some $s \in S$, i.e., $A = R + \sigma_s A$ for some finite set $R$. Moreover, for $t \in \cS$, write $A = F + \sigma_t A$ for some finite set $F$. Then $A = R + \sigma_s A = R + \sigma_s (F + \sigma_\sigma A) = R + \sigma_s F + \sigma_{st} A$. Hence it follows that $\# (A / \sigma_{st} A) < \infty$.
\eremark

In this paper, we will only consider algebraic actions $\sigma\colon S\acts A$ which satisfy \eqref{FI}. Let $\bd$ be the compact space of characters on the semilattice $\cE:=\{b+C : C\in\cC, b\in A\}\cup\{\emptyset\}$ (see \cite[\S~3.4]{BL2}). Each $\bm{x}=(x_C+C)_C\in\ol{A}$ determines an element $\chi_{\bm{x}}$ of $\bd$ by 
\[
\chi_{\bm{x}}(b+C):=\begin{cases}
1 & \text{ if } x_C+C=b+C,\\
0 & \text{ if } x_C+C\neq b+C.
\end{cases}
\]
Since \eqref{FI} is satisfied, it is not hard to see that the map $\ol{A}\to \bd$ given by $\bm{x}\mapsto \chi_{\bm{x}}$ is a homeomorphism.

In addition, we will always assume that $\sigma\colon S\acts A$ has a \emph{globalization} $\tilde{\sigma} \colon \mathscr{S} \acts \mathscr{A}$, i.e., $S$ embeds into the group $\mathscr{S}$, $A$ is a subgroup of the group $\mathscr{A}$, and $\tilde{\sigma}\colon \mathscr{S}\to\Aut(\mathscr{A})$ is an algebraic action such that $\tilde{\sigma}_s\vert_A=\sigma_s$ for all $s\in S$. Then $\mathscr{A}\rtimes \mathscr{S}$ acts on $\mathscr{A}$ by affine maps: $(z,\gamma).x=z+\tilde{\sigma}_\gamma(x)$. Reducing to $A\subseteq \mathscr{A}$, we get a partial action (in the sense of \cite{ExelPAMS}) on the group $A$. Explicitly, $g \in \mathscr{A}\rtimes \mathscr{S}$ acts by the partial bijection
\[
 A\cap g^{-1}.A\to (g.A)\cap A, \quad x\mapsto g.x.
\]
Note that $s\in S$ acts via $\sigma_s$ (where we view both maps as partial bijections on $A$). Consider the condition 
\begin{equation}
    \label{JF}\tag{JF}
    C\subseteq \ker(\id-\ti{\sigma}_g)\implies g=1\text{ for all }C\in\cC,\, g\in \gp{S},
\end{equation}
where $\gp{S}$ is the subgroup of $\msc{S}$ generated by $S$ (cf. \cite[\S~3.3]{BL2}). Moreover, our standing assumptions in this paper include that \eqref{JF} is satisfied. In this case, the partial action $\mathscr{A}\rtimes \mathscr{S} \acts A$ extends uniquely to a partial action $\mathscr{A}\rtimes \mathscr{S} \acts \ol{A}$. Given $g \in \mathscr{A}\rtimes \mathscr{S}$, let $U_{g^{-1}}\subseteq \ol{A}$ be the domain of $g$, and for $\bm{x}\in U_{g^{-1}}$, we let $g.\bm{x}$ denote the image of $\bm{x}$ under $g$ with respect to the action $\mathscr{A}\rtimes \mathscr{S} \acts \ol{A}$. The associated transformation groupoid 
\[
 (\mathscr{A}\rtimes \mathscr{S}) \ltimes\ol{A}:=\{(g,\bm{x})\in (\mathscr{A}\rtimes \mathscr{S}) \times\ol{A} :\bm{x}\in U_{g^{-1}}, \, g.\bm{x} \in \ol{A}\}
\]
is canonically isomorphic to the groupoid $\cG_{\sigma} = I_\sigma\ltimes\bd$ from \cite[\S~3]{BL2} (see \cite[\S~3.5]{BL2}). Since \eqref{FI} is satisfied, $(\mathscr{A}\rtimes \mathscr{S}) \ltimes\ol{A}$ is minimal by \cite[Corollary~7.4]{BL2}. By \cite[Theorem~4.14]{BL2}, $(\mathscr{A}\rtimes \mathscr{S}) \ltimes\ol{A}$ is effective if and only if $\sigma\colon S\acts A$ is exact in the sense of \cite[Definition~4.11]{BL2}, i.e., $\bigcap_{C\in\cC}C=\{0\}$.

\bremark
If $S$ is left Ore, then we can replace \eqref{JF} by the condition that $C \subseteq \ker(\sigma_s - \sigma_{s'})$ for some $C \in \cC$ implies that $s=s'$ (see \cite[Example~3.17~(iii)]{BL2}).
\eremark

\bremark
\label{rmk:Oregrpd}
If $S$ is left Ore, then there is a canonical partial action of $G$ on $\ol{A}$ in general, without the assumption that \eqref{JF} holds: We first construct the enveloping action as in \cite{CuLi}. $\sigma$ extends to an action of $S$ of $\ol{A}$, also denoted by $\sigma$. Then set $S^{-1} \ol{A}\defeq \ilim_S \gekl{\ol{A}, \sigma}$; this is a locally compact (non-compact) topological group. Extend $\sigma$ to $\sigma: \: \spkl{S} \curvearrowright S^{-1} \ol{A}$. This way, we obtain a global dynamical system $S^{-1}A \rtimes \spkl{S} \curvearrowright S^{-1} \ol{A}$. Then $\ol{A}$ is a clopen subset of $S^{-1} \ol{A}$, so that we obtain the desired partial dynamical system by restricting $S^{-1}A \rtimes \spkl{S} \curvearrowright S^{-1} \ol{A}$ to $\ol{A}$.

However, note that \eqref{JF} holds automatically in this setting if $A$ is torsion-free by \cite[Proposition~7.5]{BL2}.
\eremark

\bremark
\label{rem:mscA=CUPA}
We can and will always assume that $\msc{S}$ is generated by $S$, i.e., $\msc{S} = \gp{S}$. Moreover, by \cite[Proposition~2.7]{BL2}, if $S$ embeds into the group $\msc{S}$, then we can always take $\msc{A} = \Zz \msc{S} \otimes_{\Zz S} A$, and the map $A \to \Zz \msc{S} \otimes_{\Zz S} A, \, a \ma 1 \otimes a$ will always be injective if $\sigma$ admits a globalization. In this case, we have $\msc{A} = \gp{\bigcup_{s \in \msc{S}} s.(1 \otimes A)}$. Hence, for any globalization $\ti{\sigma}: \: \msc{S} \acts \msc{A}$, we may and will always assume that
\begin{equation}
\label{e:mscA=CUPA}
 \msc{A} = \gp{\textstyle{\bigcup_{s \in \msc{S}}} \ti{\sigma}_s(A)}.
\end{equation}
If $S$ is left Ore, then we can always take $\msc{S} = S^{-1} S$, and in this case, we can and will always arrange that $\msc{A} = \bigcup_{r \in S} \ti{\sigma}_r^{-1}(A)$.
\eremark

\subsection{Further properties}
\label{ss:FurtherProp}

Let us now discuss a few properties which are not part of our standing assumptions, but which we will assume for some of our results.

The principal $S$-constructible subgroups are cofinal in $\cC$ if
\begin{equation}
\label{PC} \tag{PC}
\text{for every } C\in\cC, \text{ there exists } s\in S \text{ such that } \sigma_sA\subseteq C.
\end{equation}
Note that \eqref{PC} is satisfied if $S$ is left reversible (see \cite[Proposition~7.12]{BL2}).

Consider the following condition on the algebraic action $\msc{S}\acts \msc{A}$:
\begin{equation}
\label{F}\tag{F}
\text{For all } 1\neq s\in \mathscr{S},\, 1-\ti{\sigma}_s:=\id-\ti{\sigma}_s\colon \mathscr{A} \to \mathscr{A} \text{ is injective.}
\end{equation}
Condition \eqref{F} is a freeness condition, modulo the fact that in the linear setting $0$ will always be a fixed point.

\section{From cocycles to embeddings}
\label{s:c}

Let $\sigma\colon S\acts A$ and $\tau\colon T\acts B$ be algebraic actions with globalizations $\ti{\sigma} \colon \mathscr{S} \acts \mathscr{A}$ and $\ti{\tau} \colon \msc{T} \acts \mathscr{B}$, respectively, which satisfy all our standing assumptions from \S~\ref{ss:Standing} (and we use the same notation as in \S~\ref{sec:prelim}). We will often view $S$ as a submonoid of $\msc{A} \rtimes \msc{S}$ via $S\to \msc{A} \rtimes \msc{S}$, $s\mapsto (0,s)$ and $A$ as a subgroup of $\msc{A} \rtimes \msc{S}$ via $A\to \msc{A} \rtimes \msc{S}$, $a \mapsto (a,1)$. Moreover, we will use multiplicative notation for $\msc{A} \rtimes \msc{S}$.

Suppose $\mfc \colon (\msc{A} \rtimes \msc{S})\ltimes\ol{A}\to \msc{B}\rtimes\msc{T}$ is a continuous cocycle (i.e., a groupoid homomorphism) such that $\mfc^{-1}(0,1)=\ol{A}$. Note that $\mfc$ satisfies the cocycle identity $\mfc(gh,\bm{x}) = \mfc(g,h.\bm{x}) c(h,\bm{x})$ whenever these expressions make sense. We denote the map $A\rtimes S\to \msc{B}\rtimes\msc{T}, \, p \ma \mfc(p,\bm{0})$ again by $\mfc$.

\blemma
\label{lem:constcoset}
For each $p \in A \rtimes S$, there is $C(p)\in\cC$ such that $\mfc(p,-)\colon \ol{A} \to \msc{B}\rtimes\msc{T}$ is constant on $x+\overline{C(p)}$ for every $x\in A$.
\elemma
\setlength{\parindent}{0cm} \setlength{\parskip}{0cm}

\bproof
Since $\ol{A}$ is compact and $\mfc$ is continuous, the image of $\mfc(p,-)$ must be a finite set, which we shall denote by $F$. For each $f\in F$, let $U_f:=\mfc(p,-)^{-1}(\{f\})$. Then each $U_f$ is compact open, and we have a partition $\ol{A}=\coprod_{f\in F} U_f$. Since the collection $\{x+\ol{C} : C\in\cC, x\in A\}$ forms a base consisting of compact open sets for the topology of $\ol{A}$, we can write $U_f$ as a finite disjoint union $U_f=\coprod_i x_i+\ol{C}_i$. If we now set $C_f:=\bigcap_i C_i$, then since $C_f$ is a finite index subgroup of each $C_i$, we can even write $U_f$ as a disjoint union of the form $\coprod_i y_i+\ol{C_f}$.
Now set $C(p):=\bigcap_{f\in F} C_f$, so that $\ol{A}$ is a finite disjoint union $\ol{A}=\coprod_k x_k+\overline{C(p)}$. For all $x\in A$, there exists $k$ such that $x+\overline{C(p)} = x_k+\overline{C(p)}$, and $\mfc(p,-)$ is constant on $x_k+\overline{C(p)}$.
\eproof
\setlength{\parindent}{0cm} \setlength{\parskip}{0.5cm}

\blemma
\label{lem:c_addconst_gen}
For every finitely generated subgroup $\fA$ of $A$, there exists $C_{\fA} \in\cC$ such that $\mfc(x,-)$ is constant on $a+\ol{C_{\fA}}$ for all $a\in A$ and $x\in \fA$.
\elemma
\setlength{\parindent}{0cm} \setlength{\parskip}{0cm}

\bproof
Let $a\in A$, and let $\{x_i\}_i$ be a finite collection of elements in $\fA$ that generate $\fA$ as an additive monoid. Set $C_{\fA}:=\bigcap_i C(x_i,1)$ where the subgroups $C(x_i,1)$ are provided by Lemma~\ref{lem:constcoset}. We now show that $\mfc((x,1),-)$ is constant on $a+\ol{C_{\fA}}$ by induction on the word length $\ell(x)$ of $x$ with respect to $\{x_i\}_i$. The induction base case $\ell(x)=1$ follows from Lemma~\ref{lem:constcoset}. Now suppose the claim is true for all $x\in \fA$ with $\ell(x)=l$. Given $x\in \fA$ with $\ell(x)=l+1$, we can write $x=x_ix'$ for some index $i$ and $x'\in \fA$ with $\ell(x')=l$. Now we have for all $a+\bm{x},a+\bm{y}\in a+\ol{C_{\fA}}$,
\begin{eqnarray*}
	\mfc((x,1),a+\bm{x}) &=& \mfc((x_i,1)(x',1),a+\bm{x})\\
	&=& \mfc((x_i,1),(x',1).(a+\bm{x})) \, \mfc((x',1),a+\bm{x}) \\
	&=& \mfc((x_i,1),(x',1).(a+\bm{x})) \, \mfc((x',1),a+\bm{y})\\
	&=& \mfc((x_i,1), x'.(a+\bm{x})) \, \mfc((x',1),a+\bm{y}) \\
	&=& \mfc((x_i,1), x'.(a+\bm{y})) \, \mfc((x',1),a+\bm{y}) = \mfc((x,1),a+\bm{y}).
\end{eqnarray*}
Here, we used the induction hypothesis for the third equality and Lemma~\ref{lem:constcoset} for the fifth equality ($x'.(a+\bm{x})$ and $x'.(a+\bm{y})$ both lie in $x'+a+\ol{C_{\fA}}$).
\eproof
\setlength{\parindent}{0cm} \setlength{\parskip}{0.5cm}

\bremark Lemma~\ref{lem:c_addconst_gen} can also be derived from the general results in \cite{CM}, but we chose to give a direct proof in our special setting.
\eremark

\blemma
\label{lem:fisubgroup}
For every finitely generated subgroup $\fA$ of $A$, the map $\fA \cap C_{\fA}\to \msc{B}\rtimes\msc{T}$ given by $x \ma \mfc(x,\bm{0})$ is an injective group homomorphism.
\elemma
\setlength{\parindent}{0cm} \setlength{\parskip}{0cm}

\bproof
Additivity is easy to see. Now suppose we have $x,y\in \fA \cap C_{\fA}$ with $\mfc(x,\bm{0})=\mfc(y,\bm{0})$. Consider the element $(x,\bm{0})(y,\bm{0})^{-1}$ of the groupoid $(\msc{A} \rtimes \msc{S})\ltimes \ol{A}$. Since $\mfc$ is a groupoid homomorphism, $\mfc((x,\bm{0})(y,\bm{0})^{-1})=(0,1)$. Hence, using the assumption $\mfc^{-1}(0,1)=\ol{A}$, we conclude that $(xy^{-1},\bm{0})=(x,\bm{0})(y,\bm{0})^{-1}\in\ol{A}$, which implies $x=y$. 
\eproof

Note that $\fA \cap C_{\fA}$ is a finite index subgroup of $\fA$. In particular, if $A$ is finitely generated, then we get an injective group homomorphism from a finite index subgroup of $A$ into $\msc{B}\rtimes\msc{T}$.
\setlength{\parindent}{0cm} \setlength{\parskip}{0.5cm}

\subsection{The additive homomorphism}

Given a finitely generated subgroup $\fA$ of $A$, we now want to find $\mfl \in \Zz_{>0}$ and $C \in \cC$ together with a homomorphism $\mfb \colon \fC \defeq \mfl (\fA \cap C) \to \mathscr{B}$ such that $\mfc(x) = (\mfb(x),1)$ for all $x \in \fC$.

\bprop
\label{prop:GammaAlphaEq}
Suppose that $\cA \subseteq A$ is a finite rank subgroup and that $s \in S$ satisfies $\sigma_s(\cA) \subseteq \cA$. Let $d$ be the degree of the polynomial $\det(z - \dot{\sigma}_s)$, where $\dot{\sigma}_s:=\id_{\Qz} \otimes (\sigma_s\vert_\cA)\colon \Qz \otimes \cA \to \Qz \otimes \cA$, and let $\kappa_d\in\Zz_{>0}$ be the smallest positive integer such that $p(z) \defeq \kappa_d \det(z - \dot{\sigma}_s)$ has integer coefficients, and write $p(z) = \kappa_d z^d - \kappa_{d-1} z^{d-1} - \dotso - \kappa_1 z - \kappa_0$ (for some $\kappa_{\bullet} \in \Zz$).

Then for every finitely generated subgroup $\fA$ of $\cA$, there exists $\check{C} \in \cC$ depending on $s$ such that 
\setlength{\parindent}{0cm} \setlength{\parskip}{0cm}

\begin{enumerate}[\upshape(i)]
\item The restriction of $\mfc$ to $\check{\fC} \defeq \fA \cap \check{C}$ is an injective group homomorphism $\check{\fC} \to \msc{B}\rtimes\msc{T}$. 
\item $\mfc(s)^d \mfc(x)^{\kappa_d} = \mfc(x)^{\kappa_0} \mfc(s) \dotsm \mfc(x)^{\kappa_{d-1}} \mfc(s)$ for all $x \in \check{\fC}$.
\item The following holds for all $x \in \check{\fC}$: If $\mfc(x) = (\beta, \alpha)$ and $\mfc(s) = (\delta, \gamma)$, then 
\begin{equation}
\label{e:GammaAlphaEq}
\gamma^d \alpha^{\kappa_d} = \alpha^{\kappa_0} \gamma\alpha^{\kappa_1} \gamma \dotsm \alpha^{\kappa_{d-1}} \gamma
\end{equation}
in $\mathscr{T}$.
\end{enumerate}
\eprop
\setlength{\parindent}{0cm} \setlength{\parskip}{0cm}

\bproof
Let us denote the restriction of $\sigma_s$ to $\cA$ again by $\sigma_s$. For all $x \in \cA$, the following holds in $A \rtimes S$ as $\kappa_d \sigma_s^d(x) = \kappa_{d-1} \sigma_s^{d-1}(x) + \dotso + \kappa_0 x$ in $\cA$:
$$
s^d (\kappa_d x) = \kappa_d \sigma_s^d(x) s^d = (\kappa_0 x) s (\kappa_1 x) s \dotsm (\kappa_{d-1}x) s.
$$
Given a finitely generated subgroup $\fA$ of $\cA$, choose $C_{\fA}$ as in Lemma~\ref{lem:fisubgroup}. Set
$$
 \check{C} \defeq C(s) \cap \sigma_s^{-1} C(s) \cap \dotso \cap \sigma_s^{-(d-1)} C(s) \cap C_{\fA} \cap \sigma_s^{-1} C_{\fA} \cap \dotso \cap \sigma_s^{-(d-1)} C_{\fA}.
$$
Then (i) is satisfied because of Lemma~\ref{lem:fisubgroup}.
\setlength{\parindent}{0cm} \setlength{\parskip}{0.5cm}

We have for all $x \in \check{C}$: 
\begin{align*}
\mfc(s^d (\kappa_d x)) &= \mfc(s^d (\kappa_d x), \bm{0}) = \mfc(s^d,\kappa_d x)\mfc(\kappa_d x,\bm{0})\\
&= \mfc(s^{d-1},\sigma_s(\kappa_d x)) \mfc(s,\kappa_d x) \mfc(\kappa_d x,\bm{0})\\
&= \dotso = \mfc(s, \sigma_s^{d-1}(\kappa_d x)) \mfc(s,\sigma_s^{d-2}(\kappa_d x)) \dotsm \mfc(s,\kappa_d x) \mfc(\kappa_d x,\bm{0}) \\
&= \mfc(s,\bm{0}) \mfc(s,\bm{0}) \dotsm \mfc(s,\bm{0}) \mfc(\kappa_d x,\bm{0}) = \mfc(s)^d \mfc(\kappa_d x).
\end{align*}
Here we are allowed to replace $\sigma_s^\ell(\kappa_d x)$ by $\bm{0}$ for $\ell=0,1,...,d-1$ because $x$ lies in $C(s) \cap \sigma_s^{-1} C(s) \cap \dotso \cap \sigma_s^{-(d-1)} C(s)$. 

We also have for all $x \in \fA \cap \check{C}$: 
\begin{align*}
& \ \mfc((\kappa_0 x) s (\kappa_1 x) s \dotsm (\kappa_{d-1}x) s) = \mfc((\kappa_0 x) s (\kappa_1 x) s \dotsm (\kappa_{d-1}x) s,\bm{0})\\
= & \ \mfc((\kappa_0 x) s (\kappa_1 x) s \dotsm (\kappa_{d-1}x),\bm{0}) \mfc(s,\bm{0})\\
= & \ \mfc((\kappa_0 x) s (\kappa_1 x) s \dotsm, \kappa_{d-1}x) \mfc(\kappa_{d-1}x,\bm{0}) \mfc(s,\bm{0})\\
= & \ \mfc(\kappa_0 x, \sigma_s (\kappa_1 x) + \dotso + \sigma_s^{d-1}(\kappa_{d-1}x)) \mfc(s,(\kappa_1 x) + \dotso + \sigma_s^{d-2}(\kappa_{d-1}x))\\
& \ \dotsm \\
& \ \mfc(\kappa_{d-2}x, \sigma_s(\kappa_{d-1}x)) \mfc(s, \kappa_{d-1}x)\\
& \ \mfc(\kappa_{d-1}x, \bm{0}) \mfc(s, \bm{0})\\
= & \ \mfc(\kappa_0 x, \bm{0}) \mfc(s, \bm{0}) \dotsm \mfc(\kappa_{d-2}x, \bm{0}) \mfc(s, \bm{0}) \mfc(\kappa_{d-1}x, \bm{0}) \mfc(s, \bm{0})\\
= & \ \mfc(\kappa_0 x) \mfc(s) \dotsm \mfc(\kappa_{d-1}x) \mfc(s).
\end{align*}
Here we are allowed to replace the second arguments of $\mfc$ by $\bm{0}$ because $x$ lies in $C(s) \cap \dotso \cap \sigma_s^{-(d-1)} C(s) \cap \fA \cap C_{\fA} \cap \sigma_s^{-1} C_{\fA} \cap \dotso \cap \sigma_s^{-(d-1)} C_{\fA}$.

Moreover, we have for all $x \in \check{C}$ and $\kappa \in \Zz$ that $\mfc(\kappa x) = \mfc(x)^{\kappa}$ because $x$ lies in $\fA \cap C_{\fA}$. So all in all, we have $\mfc(s)^d \mfc(x)^{\kappa_d} = \mfc(x)^{\kappa_0} \mfc(s) \dotsm \mfc(x)^{\kappa_{d-1}} \mfc(s)$ for all $x \in \check{\fC} = \fA \cap \check{C}$. This shows (ii).

Now let us prove (iii). Suppose that $\mfc(x) = (\beta,\alpha)$ and $\mfc(s) = (\delta,\gamma)$. Comparing $\mathscr{T}$-components, we obtain 
$\gamma^d \alpha^{\kappa_d} = \alpha^{\kappa_0} \gamma \alpha^{\kappa_1} \gamma \dotsm\alpha^{\kappa_{d-1}}\gamma.$
\eproof

\begin{remark}
Suppose we are in the setting of Proposition~\ref{prop:GammaAlphaEq}, and put $\epsilon \defeq \kappa_d - \kappa_{d-1} - \dotso \kappa_1 - \kappa_0$. If $\msc{T}$ is Abelian, then \eqref{e:GammaAlphaEq} is equivalent to $\alpha^\epsilon=1$. If $p(z)=z^d-\kappa_0$, then \eqref{e:GammaAlphaEq} is $\gamma^d\alpha\gamma^{-d}=\alpha^{\kappa_0}$, i.e., $\alpha$ and $\gamma^d$ satisfy the defining relation for the Baumslag--Solitar group $\textup{BS}(1,\kappa_0)\cong \Zz[1/\kappa_0]\rtimes\Zz$. If $\gamma^d$ and $\alpha$ both have infinite order, then $\gp{\gamma^d,\alpha}\cong \textup{BS}(1,\kappa_0)$.
\end{remark}
\setlength{\parindent}{0cm} \setlength{\parskip}{0.5cm}

\blemma
\label{lem:AlphaTorsion}
With the same notation as in Proposition~\ref{prop:GammaAlphaEq}, set $\epsilon \defeq \kappa_d - \kappa_{d-1} - \dotso \kappa_1 - \kappa_0$. In addition to the assumptions in Proposition~\ref{prop:GammaAlphaEq}, assume that $1 - \sigma_s$ is injective on $\cA$. Then we have $\epsilon \neq 0$. If, in addition, every $2$-generated subgroup of $\mathscr{T}$ is free or Abelian, then $\alpha^{\epsilon} = 1$.
\elemma
\setlength{\parindent}{0cm} \setlength{\parskip}{0cm}

\bproof
The first claim follows from $\kappa_d \det(1 - \dot{\sigma}_s) = \epsilon$. For the second claim, our assumption implies that the subgroup $\gp{\alpha,\gamma}\subseteq\gp{T}$ is free or Abelian. We claim that $\alpha$ and $\gamma$ must commute. Indeed, it suffices to treat the case that the subgroup is free. Because $\alpha$ and $\gamma$ satisfy the non-trivial relation \eqref{e:GammaAlphaEq}, the group $\gp{\alpha,\gamma}$ is either trivial or infinite cyclic; if it is trivial, we are done. Suppose $\gp{\alpha,\gamma}$ is cyclic. Then, in particular, $\alpha\gamma=\gamma\alpha$, as desired. Now \eqref{e:GammaAlphaEq} becomes $\gamma^d\alpha^{\kappa_d}=\gamma^d\alpha^{\kappa_0+\kappa_1+\cdots\kappa_{d-1}}$, which implies that $\alpha^{\epsilon}=1$. 
\eproof

\bex
Let us mention two classes of groups whose 2-generated subgroups are either free or Abelian. 
A group is called \emph{semifree} if it has a presentation where the only relators are of the form $[s,t]$, where $s$ and $t$ are generators. If $s,t$ are elements in a semifree group and $[s,t]\neq 1$, then $\{s,t\}$ is a basis for a free group by \cite[Theorem~1.2]{Baudisch}. 
\setlength{\parindent}{0.5cm} \setlength{\parskip}{0cm}

A group is \emph{2-free} if every subgroup generated by 2 elements is free. Given a non-empty set $\omega$ of prime numbers, a $D_\omega$-free group is a group whose elements each have exactly one $p$-th root for all primes $p\in\omega$ (see the introduction in \cite{G.Baumslag}). By \cite[\S~8]{B.Baumslag}, every $D_\omega$-free group from \cite{G.Baumslag} is 2-free.
\eex
\setlength{\parindent}{0cm} \setlength{\parskip}{0.5cm}

Let us introduce the following conditions.
\bdefin
Let $\cA \subseteq A$ be a subgroup of finite rank. Consider the following conditions:
\setlength{\parindent}{0cm} \setlength{\parskip}{0cm}

\begin{enumerate}
    \item[\upshape{(i}$_1$\upshape{)}] There exists $s \in S$ such that $\sigma_s(\cA) \subseteq \cA$, $1 - \sigma_s$ is injective on $\cA$, and every $2$-generated subgroup of $\mathscr{T}$ is free or Abelian.
    \item[\upshape{(i}$_2$\upshape{)}] There exists $s \in S$ with $\sigma_s = \kappa \, \id_A$ for some $\kappa \in \Zz \setminus \gekl{0,1}$, and for all $\alpha \in \msc{T}$, if $\alpha$ is conjugate to $\alpha^{\kappa}$ in $\msc{T}$, then $\alpha$ must be torsion.
\end{enumerate}
\edefin
\setlength{\parindent}{0cm} \setlength{\parskip}{0.5cm}

\bcor
\label{cor:lACtoB}
If $\cA \subseteq A$ is a subgroup of finite rank and (i$_1$) or (i$_2$) holds, then for all finite generated subgroups $\fA \subseteq \cA$, there exist $\mfl \in \Zz_{>0}$ and $C \in \cC$ such that, with $\fC \defeq \mfl (\fA \cap C)$, there exists an injective homomorphism $\mfb: \: \fC \to \msc{B}$ such that $\mfc(x) = (\mfb(x),1)$ for all $x \in \fC$.
\ecor
\setlength{\parindent}{0cm} \setlength{\parskip}{0cm}

\bproof
Let $s \in S$ be as in (i$_1$) or (i$_2$), and $\check{C}$, $\check{\fC}$ as in Proposition~\ref{prop:GammaAlphaEq}. Let $\mfc_{\mathscr{T}}$ be the composition $A \to \msc{B}\rtimes\msc{T} \onto \mathscr{T}$, where the first map is $\mfc$ and the second map is the canonical projection. Now $\mfc_{\mathscr{T}}(\check{\fC})$ is finitely generated. Moreover, if (i$_1$) holds, then Lemma~\ref{lem:AlphaTorsion} implies that $\mfc_{\mathscr{T}}(\check{\fC})$ is torsion, hence finite. If (i$_2$) holds, then Proposition~\ref{prop:GammaAlphaEq}~(iii) implies that $\mfc_{\mathscr{T}}(\check{\fC})$ is torsion, hence finite. Set $\mfl \defeq \# \mfc_{\mathscr{T}}(\check{\fC})$. Then for all $x \in \mfl \check{\fC}$, $\mfc_{\mathscr{T}}(x) = 1$, so that our claim follows (Lemma~\ref{lem:fisubgroup} gives injectivity of $\mfb$).
\eproof
\setlength{\parindent}{0cm} \setlength{\parskip}{0.5cm}

Recall that $\ti{\tau} \colon \msc{T} \acts \mathscr{B}$ is our globalization of $\tau\colon S\acts B$. 

\bdefin
Let $\cA \subseteq A$ be a subgroup of finite rank. Consider the following conditions:
\setlength{\parindent}{0cm} \setlength{\parskip}{0cm}

\begin{enumerate}
    \item[\upshape{(ii}$_1$\upshape{)}] For all torsion orders $l>1$ of elements of $\mathscr{T}$, there exists $1 \neq s_l \in S$ and coprime integers $\mu,\nu\in\Zz_{>0}$ such that $\sigma_{s_l}(\cA) \subseteq \cA$ and $\mu \det(z - \dot{\sigma}_{s_l}) = \mu z^{\delta} - \nu$ with $\gcd(l,\mu) > 1$ or $\gcd(l,\nu) > 1$. 
    \item[\upshape{(ii}$_2$\upshape{)}] Condition \eqref{F} holds for $\ti{\tau}$.
\end{enumerate}
\edefin
Note that, in particular, (ii$_1$) holds if $\msc{T}$ is torsion-free or if for all $\kappa \in \Zz_{>0}$ there exists $s_{\kappa} \in S$ with $\sigma_{s_{\kappa}} = \kappa \, \id_A$.

\bdefin
\label{def:(III)}
We say that condition \hypertarget{(III)}{(III)} holds if there exists $s \in S$ such that $\sigma_s = \kappa \, \id_A$ for some $\kappa \in \Zz \setminus \gekl{0,1}$, $\msc{B}$ is of finite rank, and condition \eqref{F} holds for $\ti{\tau}$.
\edefin
\setlength{\parindent}{0cm} \setlength{\parskip}{0.5cm}

\bcor
\label{cor:ACtoB}
Let $\cA \subseteq A$ be a subgroup of finite rank. Assume that one of the following is true:
\setlength{\parindent}{0cm} \setlength{\parskip}{0cm}

\begin{itemize}
    \item (i$_1$) or (i$_2$) holds, and (ii$_1$) is satisfied or (ii$_2$) holds and $A$ is torsion-free,
    \item \hyperlink{(III)}{(III)} holds and $A$ is torsion-free.
\end{itemize}
Then for all finitely generated subgroup $\fA \subseteq \cA$, there exists $C \in \cC$ such that, with $\fC \defeq \fA \cap C$, there exists an injective homomorphism $\mfb: \: \fC \to \msc{B}$ such that $\mfc(x) = (\mfb(x),1)$ for all $x \in \fC$.
\ecor
\setlength{\parindent}{0cm} \setlength{\parskip}{0cm}

\bproof
To prove the first item, let $s \in S$ be as in (i$_1$) or (i$_2$), and $\check{C}$, $\check{\fC}$ as in Proposition~\ref{prop:GammaAlphaEq}. First assume that (ii$_1$) holds. Let $\mfc_{\mathscr{T}}$ be the composition $A \to \msc{B}\rtimes\msc{T} \onto \mathscr{T}$, where the first map is $\mfc$ and second map is the canonical projection. As $\mfc_{\mathscr{T}}(\check{\fC})$ is finitely generated, the set $\Lz$ of possible non-trivial torsion orders of elements in $\mfc_{\mathscr{T}}(\check{\fC})$ is finite. For each $l \in \Lz$, choose $s_l$ as in (ii$_1$) and set $C \defeq \check{C} \cap \bigcap_{l \in \cL} \check{C}(s_l)$. Applying Proposition~\ref{prop:GammaAlphaEq} to $s_l$, we obtain that for all $x \in \fA \cap C$ with $\mfc(x) = (\beta,\alpha)$, we have $\alpha^{\mu}$ and $\alpha^{\nu}$ are conjugate in $\mathscr{T}$, where $\mu$ and $\nu$ are the natural numbers from (ii$_1$). If $\alpha \neq 1$, then the torsion order of $\alpha$ is a number $l \in \Lz$, but (ii$_1$) implies that $\alpha^{\mu}$ and $\alpha^{\nu}$ have different torsion orders, which is absurd. Now suppose that (ii$_2$) holds. Take $x \in \check{\fC}$ and write $\mfc(x) = (\beta,\alpha)$. Lemma~\ref{lem:AlphaTorsion} implies that $\alpha^{\epsilon} = 1$. Then
$$
(\beta,\alpha)^{\epsilon} = (\beta + \ti{\tau}_\alpha(\beta) + \dotso + \ti{\tau}_\alpha^{\epsilon-1}(\beta),\alpha^{\epsilon}) = (\beta + \ti{\tau}_\alpha(\beta) + \dotso + \ti{\tau}_\alpha^{\epsilon-1}(\beta),1),
$$
and
$$
(1 - \ti{\tau}_\alpha)(\beta + \ti{\tau}_\alpha(\beta) + \dotso + \ti{\tau}_\alpha^{\epsilon-1}(\beta)) = 0,
$$
which implies $\beta + \ti{\tau}_\alpha(\beta) + \dotso + \ti{\tau}_\alpha^{\epsilon-1}(\beta) = 0$ if $\alpha \neq 1$ by \eqref{F}. So $\mfc(x)^{\epsilon} = 1$ and hence $\mfc(\epsilon x) = 1$ and thus $\epsilon x = 0$. As $A$ is torsion-free, this implies $x = 0$. So if $x \neq 0$, we must have $\alpha = 1$.
\setlength{\parindent}{0cm} \setlength{\parskip}{0.5cm}

Now we prove the second item. As above, let $s \in S$ be as in (III), and $\check{C}$, $\check{\fC}$ as in Proposition~\ref{prop:GammaAlphaEq}. Given $x \in \check{\fC}$, write $\mfc(x) = (\beta,\alpha)$ and $\mfc(s) = (\delta,\gamma)$. Equation \eqref{e:GammaAlphaEq} implies that $\gamma \alpha = \alpha^{\kappa} \gamma$, and therefore $\alpha = \gamma^{-1} \alpha^{\kappa} \gamma$. It follows that every eigenvalue of $\dot{\tau}_\alpha \defeq \id_\Qz\otimes \tau_\alpha$ is a root of unity. Indeed, take $\lambda_1 \in \Sp(\dot{\tau}_\alpha)$. Then $\alpha = \gamma^{-1} \alpha^{\kappa} \gamma$ implies that there exists $\lambda_2 \in \Sp(\dot{\tau}_\alpha)$ with $\lambda_1 = \lambda_2^{\kappa}$. Similarly, there exist $\lambda_3, \lambda_4, \dotsc \in \Sp(\dot{\tau}_\alpha)$ such that $\lambda_i = \lambda_{i+1}^{\kappa}$. As $\Sp(\dot{\tau}_\alpha)$ is finite, we must have $\lambda_i = \lambda_{i+p}$ for some $i$ and $p$. It follows that $\lambda_i = \lambda_i^{\kappa^p}$ and hence that $\lambda_i$ is a root of unity. But then $\lambda_1 = \lambda_i^{\kappa^i}$ implies that $\lambda_1$ must be a root of unity as well. Hence there exists $m\in\Zz_{>0}$ such that $1$ is an eigenvalue of $\dot{\tau}_\alpha^m$. This implies that $1-\dot{\tau}_\alpha^m$ is not injective, so that $\dot{\tau}_\alpha^m=1$ by \eqref{F}. Now argue as for the first item that this -- together with \eqref{F} -- implies $\alpha = 1$.
\eproof
\setlength{\parindent}{0cm} \setlength{\parskip}{0.5cm}

\bremark
The difference between Corollaries~\ref{cor:lACtoB} and \ref{cor:ACtoB} is that in the latter, we may choose $\mfl = 1$.
\eremark

\bdefin
\label{def:I,II}
Assume that $A = \bigcup_n \cA_n$ for an increasing family of finite rank subgroups $\cA_n$. Consider the following conditions:
\setlength{\parindent}{0cm} \setlength{\parskip}{0cm}

\begin{enumerate}
    \item[\upshape{\hypertarget{(I)}{(I)}}] Condition (i$_1$) holds for $\cA_n$ for all $n$, or condition (i$_2$) holds.
    \item[\upshape{\hypertarget{(II)}{(II)}}] Condition (ii$_1$) holds for $\cA_n$ for all $n$, or condition (ii$_2$) holds and $A$ is torsion-free.
\end{enumerate}
\edefin
\setlength{\parindent}{0cm} \setlength{\parskip}{0.5cm}

\bcor
\label{cor:increasingfam}
Suppose we can write $A=\bigcup_n \cA_n$ for an increasing family $\cA_n\subseteq A$ of finite rank subgroups. 

If \hyperlink{(I)}{(I)} holds, then 
\setlength{\parindent}{0cm} \setlength{\parskip}{0cm}

\begin{enumerate}
    \item[\upshape{(a)}] there is an increasing family of finitely generated subgroups $\fA_k\subseteq A$ with $A=\bigcup_k\fA_k$, and for any such $\fA_k$, there are $\mfl_k \in \Zz_{>0}$ and $C_k \in \cC$ such that, with $\fC_k \defeq \mfl_k (\fA_k \cap C_k)$, there exists an injective homomorphism $\mfb_k: \: \fC_k \to \msc{B}$ such that $\mfc(x) = (\mfb_k(x),1)$ for all $x \in \fC_k$.
\end{enumerate}
\setlength{\parindent}{0cm} \setlength{\parskip}{0.5cm}

If \hyperlink{(I)}{(I)} and \hyperlink{(II)}{(II)} hold, or if \hyperlink{(III)}{(III)} is satisfied and $A$ is torsion-free, then
\setlength{\parindent}{0cm} \setlength{\parskip}{0cm}

\begin{enumerate}
    \item[\upshape{(a*)}] there is an increasing family of finitely generated subgroups $\fA_k\subseteq A$ with $A=\bigcup_k\fA_k$, and for any such $\fA_k$, there is $C_k \in \cC$ such that, with $\fC_k \defeq \fA_k \cap C_k$, there exists an injective homomorphism $\mfb_k: \: \fC_k \to \msc{B}$ such that $\mfc(x) = (\mfb_k(x),1)$ for all $x \in \fC_k$.
\end{enumerate}
\ecor
\setlength{\parindent}{0cm} \setlength{\parskip}{0cm}

\bproof
Since $A=\bigcup_n \cA_n$, we can find $\fA_k$ with the desired properties. Moreover, for each $k$, since $\fA_k$ is finitely generated, and $A=\bigcup_n \cA_n$, we can find $n$ such that $\fA_k\subseteq\cA_n$. Now apply Corollaries~\ref{cor:lACtoB} and \ref{cor:ACtoB}.
\eproof

\bremark
\label{rmk:index}
With the notation from Corollary~\ref{cor:increasingfam}, we have $\fA_k/\fC_k\hookrightarrow A/C_k$, so that $\#(\fA_k/\fC_k)$ is finite and divides $\#(A/C_k)$. 
\eremark
\setlength{\parindent}{0cm} \setlength{\parskip}{0.5cm}

\subsection{The multiplicative homomorphism}

Define $\mft \colon \mathscr{S} \to \mathscr{T}$ to be the composition $\mathscr{S} \to \mathscr{B} \rtimes \mathscr{T} \onto \mathscr{T}$, where the first arrow is given by $\mfc(-,\bm{0})$, and the second arrow is the canonical projection. Clearly, $\mft$ is a homomorphism.

\blemma
\label{lem:alpha-inj}
\begin{enumerate}[\upshape(i)]
\item If $A$ is torsion-free and (a) holds, then $\mft$ is injective.
\item Suppose that condition \eqref{F} is satisfied for $\ti{\sigma}$. If there exist $\mfl \in \Zz_{>0}$, a non-zero, finitely generated subgroup $\fA \subseteq A$, $C \in \cC$ and an injective homomorphism $\mfb: \: \fC \defeq \mfl(\fA \cap C) \to \msc{B}$ such that $\mfc(x) = (\mfb(x),1)$ for all $x \in \fC$, then $\mft$ is injective.
\end{enumerate}
\elemma
\setlength{\parindent}{0cm} \setlength{\parskip}{0cm}

\bproof
Suppose $\mft(s) = \mft(s')$. Then $\mfc(s^{-1}s', \bm{0}) = (\chi,\mft(s^{-1}s')) = (\chi,1)$ for some $\chi \in \mathscr{B}$. Set $\varepsilon \defeq s^{-1}s'$. Since $0$ lies in the domain of $\ti{\sigma}_\varepsilon$, and since $\mfc(\varepsilon,-)$ is locally constant, there must exist $C(\varepsilon) \in \cC$ such that $C(\varepsilon)$ is contained in the domain of $\ti{\sigma}_\varepsilon$ and $\mfc(\varepsilon,-)$ is constant on $\overline{C(\varepsilon)}$. Now suppose that $\fC \subseteq A$ is a non-zero subgroup such that there exists an injective homomorphism $\mfb: \: \fC \to \msc{B}$ with $\mfc(x) = (\mfb(x),1)$ for all $x \in \fC$. Then $\mfc(\varepsilon,\bm{0})$ commutes with $\mfc(x,\bm{0})$ for all $x \in \fC$. We have $(0,\varepsilon)(x,1) = (\ti{\sigma}_\epsilon(x),\epsilon) = (\ti{\sigma}_\varepsilon(x),1)(0,\varepsilon)$. So if $0 \neq x \in \fC \cap C(\varepsilon)$, then 
$$
 \mfc((\ti{\sigma}_\varepsilon(x),\varepsilon),\bm{0}) = \mfc(\ti{\sigma}_\varepsilon(x) \varepsilon,\bm{0}) = \mfc(\ti{\sigma}_\varepsilon(x),\bm{0}) \mfc(\varepsilon,\bm{0}).
$$
At the same time,
$$
 \mfc((\ti{\sigma}_\varepsilon(x),\epsilon),\bm{0}) = \mfc(\varepsilon x,\bm{0}) = \mfc(\varepsilon,x) \mfc(x,\bm{0}) = \mfc(\varepsilon,\bm{0}) \mfc(x,\bm{0}) = \mfc(x,\bm{0}) \mfc(\varepsilon,\bm{0}).
$$
Here we are allowed to replace $x$ by $\bm{0}$ because $x$ lies in $C(\varepsilon)$. Moreover, we used that $c(\varepsilon,\bm{0})$ commutes with $c(x,\bm{0})$. Therefore, a comparison yields $\mfc(\ti{\sigma}_\varepsilon(x),\bm{0}) = \mfc(x,\bm{0})$. It follows that $\ti{\sigma}_\varepsilon(x) = x$.
\setlength{\parindent}{0cm} \setlength{\parskip}{0.5cm}

For (i), write $A = \bigcup_k \fA_k$ as in (a). Applying the above to $\fC = \fC_k$ from (a), we obtain $\ti{\sigma}_\varepsilon(x) = x$ for all $x \in \fC_k \cap C(\varepsilon)$ for all $k$. Now let $y\in A$; then there exists $k$ with $y\in \fA_k$. Since $\fC_k \cap C(\varepsilon)$ is finite index in $\fA_k$, there exists $N\in\Zz_{>0}$ such that $Ny\in \fC_k \cap C(\varepsilon)$. We have $N \ti{\sigma}_\varepsilon(y)=\ti{\sigma}_\varepsilon(Ny)=Ny$, so that, because $A$ is torsion-free, $\ti{\sigma}_\varepsilon(y)=y$. Now \eqref{JF} implies $\varepsilon = 1$.

For (ii), if $\ti{\sigma}_\varepsilon(x) = x$ for any $x \neq 0$, then condition \eqref{F} implies that $\varepsilon = 1$.
\eproof

\setlength{\parindent}{0cm} \setlength{\parskip}{0.5cm}

\subsection{Equivariance}
\label{ss:Equiv}

\bdefin
We set $\Iz \defeq \menge{\#(A/C)}{C \in \cC}$.
\edefin

\blemma
\label{lem:equi}
Assume that $\scB$ is torsion-free. If (a) holds, then for all $s \in S$, $x \in \fC_k$, $\sigma_s(x) \in \fA_l$ ($k\leq l$) there exists $n \in \Zz_{>0}$ such that $n \sigma_s(x) = \sigma_s(nx) \in \fC_l$, and we have $\mfb(\sigma_s(nx)) = \ti{\tau}_{\mft(s)}(\mfb(nx))$. If (a*) holds, then we may take $n \in \Iz$ in the statement above.
\elemma
\setlength{\parindent}{0cm} \setlength{\parskip}{0cm}

\bproof
Since $\fC_l$ is of finite index in $\fA_l$ (see Remark~\ref{rmk:index}), there exists $n \in \Zz_{>0}$ such that $n \sigma_s(x) = \sigma_s(nx) \in \fC_l$ (since $l$ depends on $s$, $n$ also depends on $s$). Since $\fC_k \cap C(s)$ is of finite index in $\fC_k$, we can find $N$ such that $Nnx \in \fC_k \cap C(s)$. Set $y \defeq Nnx$. Let $\mfc(s,\bm{0}) = (\delta,\mft(s))$ for some $\delta \in \scB$. We have 
\begin{align*}
& \ (\mfb(\sigma_s(y)) + \delta,\mft(s)) = (\mfb(\sigma_s(y)),1) (\delta,\mft(s)) = \mfc(\sigma_s(y),\bm{0}) \mfc(s,\bm{0})\\
= & \ \mfc((\sigma_s(y),s),\bm{0}) = \mfc((0,s)(y,0),\bm{0}) = \mfc((0,s),y)c((y,1),\bm{0})\\
= & \ \mfc((0,s),\bm{0}) \mfc((y,1),\bm{0}) = (\delta,\mft(s)) (\mfb(y),1) = (\delta + \ti{\tau}_{\mft(s)}(\mfb(y)), \mft(s)),
\end{align*}
where the sixth equality uses that $y \in C(s)$. Hence $\mfb(\sigma_s(y)) = \ti{\tau}_{\mft(s)}(\mfb(y))$, i.e., $N \mfb(\sigma_s(nx)) = N \ti{\tau}_{\mft(s)}(\mfb(nx))$. We obtain $\mfb(\sigma_s(nx)) = \ti{\tau}_{\mft(s)}(\mfb(nx))$ as $\scB$ is torsion-free, as desired.
\eproof
\setlength{\parindent}{0cm} \setlength{\parskip}{0.5cm}

\subsection{Conclusion: The embedding theorem}
Assume that $\sigma \colon S \acts A$ and $\tau \colon T \acts B$ are algebraic actions satisfying our standing assumptions from \S~\ref{sec:prelim}. Let $\Zz[\Iz^{-1}]$ be the subring of $\Qz$ generated by $\Zz$ together with $\menge{\frac{1}{n}}{n \in \Iz}$. We start with the following observation.
\blemma
\label{lem:QA=QA}
We have $\Zz[\Iz^{-1}] \otimes A = \Zz[\Iz^{-1}] \otimes \msc{A}$, i.e., the canonical map $\Zz[\Iz^{-1}] \otimes A \to \Zz[\Iz^{-1}] \otimes \msc{A}, \, 1 \otimes x \ma 1 \otimes x$ is an isomorphism.
\elemma
\setlength{\parindent}{0cm} \setlength{\parskip}{0cm}

Recall that we always assume \eqref{e:mscA=CUPA}, i.e., $\msc{A} = \gp{\bigcup_{s \in \msc{S}} \ti{\sigma}_s(A)}$. Otherwise, Lemma~\ref{lem:QA=QA} would not be true in general.
\bproof
Because of \eqref{e:mscA=CUPA}, it suffices to prove that for all $x \in \ti{\sigma}_{t_1}^{-1} \ti{\sigma}_{s_1} \dotso \ti{\sigma}_{t_l}^{-1} \ti{\sigma}_{s_l} A$ (where $t_1, s_1, \dotsc, t_l, s_l \in S$) there exists $N$ in the multiplicative submonoid $\spkl{\Iz}^+$ of $\Zz\reg$ generated by $\Iz$ with $Nx \in A$. We proceed inductively on $l$. For the case $l=1$, observe that $[\ti{\sigma}_t^{-1} A: A] < \infty$ for all $t \in S$ because $\ti{\sigma}_t$ induces a bijection $\ti{\sigma}_t^{-1} A / A \cong A / \sigma_t A$. Now suppose that $x = \ti{\sigma}_t^{-1} \ti{\sigma}_s (y)$ for some $y \in \scA$ with $Ny \in A$ for some $N \in \spkl{\Iz}^+$. Since $[\ti{\sigma}_t^{-1} A: A] < \infty$, there exists $M \in \spkl{\Iz}^+$ with $M \ti{\sigma}_t^{-1} (\sigma_s (Ny)) \in A$. Hence $MN x \in A$, as desired.
\eproof
\setlength{\parindent}{0cm} \setlength{\parskip}{0.5cm}

\bdefin
\label{def:STAR}
We say that condition \hypertarget{(*)}{($*$)} is satisfied if $A = \bigcup_n \cA_n$ for an increasing family of finite rank subgroups $\cA_n$, conditions \hyperlink{(I)}{(I)} and \hyperlink{(II)}{(II)} hold, or condition \hyperlink{(III)}{(III)} holds; in addition, there exists an increasing family of finitely generated subgroups $\fA_k \subseteq A$ with $A = \bigcup_k \fA_k$ such that $\sigma_s(\fA_k) \subseteq \fA_k$ for all $s \in S$ and all $k$, $S$ and $T$ are right reversible, $\msc{A} = S^{-1} A$, $\msc{S} = S^{-1} S$, $\msc{B} = T^{-1} B$, $\msc{T} = T^{-1} T$, and $\sigma\colon S\acts A$ satisfies \eqref{PC}.
\edefin
Now suppose $\mfc \colon (\mathscr{A} \rtimes \mathscr{S}) \ltimes \ol{A} \to \mathscr{B} \rtimes \mathscr{T}$ is a continuous cocycle such that $\mfc^{-1}(0,1)=\ol{A}$.

\btheo
\label{thm:mainthm}
Suppose $\msc{A}$ and $\msc{B}$ are torsion-free.
\setlength{\parindent}{0cm} \setlength{\parskip}{0cm}

\begin{enumerate}[\upshape(i)]
	\item If \hyperlink{(I)}{(I)} holds, then there exist injective homomorphisms $\mft: \: \msc{S} \into \msc{T}$ and $\dot{\mfb} \colon \Qz \otimes \scA \into \Qz \otimes \scB$ such that $\dot{\mfb}(\dot{\sigma}_s(x)) = \dot{\tau}_{\mft(s)}(\dot{\mfb}(x))$ for all $s\in\msc{S}$ and $x \in \msc{A}$.
	\item If \hyperlink{(*)}{($*$)} holds, then there exist injective homomorphisms $\mft: \: \msc{S} \into \msc{T}$ and $\mfb' \colon S^{-1}A\to T^{-1}B$ such that $\mfb'(\sigma_s(x))=\ti{\tau}_{\mft(s)}(\mfb'(x))$ for all $x\in S^{-1}A$ and $s\in \scS$.
\end{enumerate}
\etheo
\setlength{\parindent}{0cm} \setlength{\parskip}{0cm}

\bproof
(i): We proceed in several steps. First, we claim that for each $k$, $\mfb_k\colon \fC_k\to\msc{B}$ has a unique extension to an injective homomorphism $\ti{\mfb}_k\colon\fA_k\to \Qz \otimes \msc{B}$, and that moreover, $\ti{\mfb}_k$ satisfies 
\begin{equation}
    \label{eqn:ext}
    \ti{\mfb}_k(x)=m^{-1}\mfb_k(mx),
\end{equation}
for any $m \in\Zz_{>0}$ with $m\fA_k\subseteq \fC_k$. To see this, choose $m\in\Zz_{>0}$ such that $m\fA_k\subseteq \fC_k$, and define $\tilde{\mfb}_k(x):=m^{-1}\mfb_k(mx)\in\Qz\otimes\msc{B}$. It is easy to check that $\tilde{\mfb}_k$ is a group homomorphism (i.e., additive) and that $\tilde{\mfb}_k$ is injective (here we need that $\fA_k$ is torsion-free). This is independent of the choice of $m$, because if $m'\in\Zz_{>0}$ with $m'\fA_k \subseteq\fC_k$, then for $x\in\fA_k$, we have
\[
m^{-1} \mfb_k(mx) = m^{-1} (m')^{-1} \mfb_k(m'mx) = (m')^{-1} m^{-1} m \mfb_k(m'x) = (m')^{-1} \mfb_k(m'x).
\]
Next, we claim that the maps $\ti{\mfb}_k$ are compatible in the sense that $\ti{\mfb}_l\vert_{\fA_k\cap\fA_l}=\ti{\mfb}_k\vert_{\fA_k\cap\fA_l}$ for all $k,l\in\Zz_{>0}$. Choose $m$ large enough so that $m\fA_k\subseteq \fC_k$ and $m\fA_l\subseteq \fC_l$. Then, using \eqref{eqn:ext}, we have for all $x\in \fA_k\cap\fA_l$, that $\ti{\mfb}_l(x)=m^{-1}\mfb_l(mx)=m^{-1}\mfb_k(mx)=\ti{\mfb}_k(x)$. It follows that we get a well-defined injective homomorphism $\ti{\mfb}\colon A=\bigcup_k\fA_k\to  \Qz \otimes \msc{B}$ such that $\ti{\mfb}\vert_{\fA_k}=\ti{\mfb}_k$ for all $k$.
\setlength{\parindent}{0cm} \setlength{\parskip}{0.5cm}

Let us show that $\ti{\mfb}$ is equivariant. Let $x\in A$ and $s\in S$. Then $x\in \fA_k$ for some $k$, so by Lemma~\ref{lem:equi}, we can find $n\in\Zz_{>0}$ large enough so that $nx\in\fC_k$, $n\sigma_s(x)=\sigma_s(nx)\in\fC_l$, and $\mfb(\sigma_s(nx))= \ti{\tau}_{\mft(s)}(\mfb(nx))$. Now we have
\[
\tilde{\mfb}(\sigma_s(x))=n^{-1}\mfb(\sigma_s(nx))=n^{-1}\ti{\tau}_{\mft(s)}(\mfb(nx))=\dot{\tau}_{\mft(s)}(n^{-1}\mfb(nx))=\dot{\tau}_{\mft(s)}(\tilde{\mfb}(x)).
\]
Lastly, we extend $\tilde{\mfb}$ to $\Qz \otimes A$ as follows: Given $\frac{p}{q} \otimes x \in \Qz \otimes A$ for $p \in \Zz$ and $q \in \Zz\reg$, there exists a unique element $y \in \msc{B}$ such that $qy = \ti{\mfb}(px)$ in $\msc{B}$, and we set $\dot{\mfb}(\frac{p}{q} \otimes x) \defeq y$. Now it is straightforward to check that this is independent of $p$, $q$ and $x$, and that $\dot{\mfb}$ is an injective homomorphism $\Qz \otimes \scA \into \Qz \otimes \scB$. Moreover, equivariance of $\ti{\mfb}$ with respect to $\ti{\sigma}$ and $\dot{\tau}$ implies equivariance of $\dot{\mfb}$ with respect to $\dot{\sigma}$ and $\dot{\tau}$.

(ii): First of all, every element of $S^{-1} A$ is of the form $s^{-1} a = \ti{\sigma}_s^{-1} a$ for some $a \in \fC_k$ (for some $k$). Indeed, the statement is clear if we just ask for $a \in \fA_k$ for some $k$ because of Remark~\ref{rem:mscA=CUPA}. By \eqref{PC}, there exists $\dot{s} \in S$ such that $\sigma_{\dot{s}}(a) \in C_k$, so that $\sigma_{\dot{s}}(a) \in \fC_k$ because $\sigma_{\dot{s}}(\fA_k) \subseteq \fA_k$, and we have $\ti{\sigma}_s^{-1} a = \ti{\sigma}_s^{-1} \ti{\sigma}_{\dot{s}}^{-1} (\ti{\sigma}_{\dot{s}} a)$, as desired.

Now given an element of $S^{-1} A$ of the form $\ti{\sigma}_s^{-1} a$ for some $a \in \fC_k$ (for some $k$), we claim that $\mfb'(\ti{\sigma}_s^{-1} a) \defeq \ti{\tau}_{\mft(s)}^{-1}(\mfb(a))$ is well defined and has the desired properties. To prove that it is well-defined, assume that $\ti{\sigma}_s^{-1} a = \ti{\sigma}_t^{-1} \dot{a}$. Since $S$ is right reversible, there exist $u, v \in S$ with $us = vt$. Hence $\ti{\sigma}_s^{-1} \ti{\sigma}_u^{-1} \ti{\sigma}_u a = \ti{\sigma}_s^{-1} a = \ti{\sigma}_t^{-1} \dot{a} = \ti{\sigma}_t^{-1} \ti{\sigma}_v^{-1} \ti{\sigma}_v \dot{a}$. It follows that $\ti{\sigma}_u a = \ti{\sigma}_v \dot{a}$. Now choose an integer $m$ such that $m \sigma_u(a), m \sigma_v(a) \in \fC_k$ and that equivariance holds (here, we use Lemma~\ref{lem:equi}). Then we have
$m \ti{\tau}_{\mft(s)}^{-1} \mfb(a) = \ti{\tau}_{\mft(s)}^{-1} \mfb(ma) = \ti{\tau}_{\mft(s)}^{-1} \ti{\tau}_{\mft(u)}^{-1} \mfb(\sigma_u(ma)) = \ti{\tau}_{\mft(t)}^{-1} \ti{\tau}_{\mft(v)}^{-1} \mfb(\sigma_v(m \dot{a})) = \ti{\tau}_{\mft(t)}^{-1} \mfb(m \dot{a}) = m \ti{\tau}_{\mft(t)}^{-1} \mfb(\dot{a})$. As $\scB$ is torsion-free, we conclude that $\ti{\tau}_{\mft(s)}^{-1} \mfb(a) = \ti{\tau}_{\mft(t)}^{-1} \mfb(\dot{a})$, as desired. It is easy to see that $\mfb'$ is additive.

To show equivariance, take $r, s \in S$. Since $S$ is right reversible, there exist $u, v \in S$ with $ur = vs$ and thus $r s^{-1} = u^{-1} v$. Given $a \in \fC_k$, choose an integer $m$ such that $\ti{\sigma}_v (ma) \in \fC_k$ and equivariance holds (here, we use Lemma~\ref{lem:equi}). Then we have
$m \mfb'(\ti{\sigma}_r \ti{\sigma}_s^{-1} (a)) = \mfb'(\ti{\sigma}_u^{-1} \ti{\sigma}_v (ma)) = \ti{\tau}_{\mft(u)}^{-1} \mfb(\ti{\sigma}_v (ma))
 = \ti{\tau}_{\mft(u)}^{-1} \ti{\tau}_{\mft(v)} \mfb(ma) = m \ti{\tau}_{\mft(r)} \ti{\tau}_{\mft(s)}^{-1} \mfb(a)
$.
As $\scB$ is torsion-free, we deduce $\mfb'(\ti{\sigma}_r \ti{\sigma}_s^{-1} (a)) = \ti{\tau}_{\mft(r)} (\ti{\tau}_{\mft(s)}^{-1} \mfb(a)) = \ti{\tau}_{\mft(r)} (\mfb'(\ti{\sigma}_s^{-1} a))$, as desired.
\eproof

\bremark
Actually we obtain an equivariant embedding $\Zz[\Iz^{-1}] \otimes \scA \into \Qz \otimes \scB$. And if \hyperlink{(I)}{(I)} and \hyperlink{(II)}{(II)} are satisfied, or \hyperlink{(III)}{(III)} holds, then we even obtain an embedding $\Zz[\Iz^{-1}] \otimes \scA \into \Zz[\Iz^{-1}] \otimes \scB$, using Lemma~\ref{lem:equi}, because (a*) holds by Corollary~\ref{cor:increasingfam}.
\eremark
\setlength{\parindent}{0cm} \setlength{\parskip}{0.5cm}

\subsection{Consequences}
\label{ss:Cons}

Let us formulate symmetrized versions of our rigidity results.

Schmidt defined \emph{finite (algebraic) equivalence} for algebraic actions of a fixed group (see, e.g., \cite[Definition~8.1]{SchPIMS}). The dual version of Schmidt's notion will appear naturally in our setting. In order to explain this, let us introduce some terminology. 

\bdefin
\begin{enumerate}[\upshape(i)]
    \item An \emph{(algebraic) embedding} of $\ti{\sigma}\colon\msc{S}\acts \msc{A}$ into $\ti{\tau}\colon\msc{T} \acts \msc{B}$ consists of a pair $(\mft,\mfb)$, where $\mft\colon \msc{S}\hookrightarrow\msc{T}$ and $\mfb\colon \msc{A}\hookrightarrow \msc{B}$ are injective homomorphisms such that $\mfb(\ti{\sigma}_s(x))=\ti{\tau}_{\mft(s)}(\mfb(x))$ for all $s\in\msc{S}$ and $x\in \msc{A}$. The embedding $(\mft,\mfb)$ is called \emph{finite index} if the image of $\mfb$ has finite index in $\msc{B}$, \emph{full} if $\mfb$ is surjective, and \emph{strict} if $\mft$ is an isomorphism.
    
    \item We say that $\ti{\sigma}\colon\msc{S}\acts \msc{A}$ and $\ti{\tau}\colon\msc{T} \acts \msc{B}$ are \emph{mutually embeddable} (written $\ti{\sigma}\colon\msc{S}\acts \msc{A} \sim_{\cM \cE} \ti{\tau}\colon\msc{T} \acts \msc{B}$) if each can be embedded into the other. If, in addition, the embeddings can be chosen to be of finite index, then we write $\ti{\sigma}\colon\msc{S}\acts \msc{A} \sim_{\cM \cE_{FI}} \ti{\tau}\colon\msc{T} \acts \msc{B}$. Given an Abelian group $\cQ$, we write $\ti{\sigma}\colon\msc{S}\acts \msc{A} \sim_{\cM \cE_{\cQ}} \ti{\tau}\colon\msc{T} \acts \msc{B}$ if $\dot{\sigma}\colon\msc{S}\acts \cQ \otimes \msc{A} \sim_{\cM \cE} \dot{\tau}\colon\msc{T} \acts \cQ \otimes \msc{B}$, where $\dot{\sigma}_s = \id_{\cQ} \otimes \ti{\sigma}_s$ and $\dot{\tau}:=\id_{\cQ} \otimes \ti{\tau}_s$.  
    We write $\ti{\sigma}\colon\msc{S}\acts \msc{A} \sim_{\cM \cE_{\cQ^{\cong}}} \ti{\tau}\colon\msc{T} \acts \msc{B}$ if there exist full embeddings of $\dot{\sigma}\colon\msc{S}\acts \cQ \otimes \msc{A}$ into $\dot{\tau}\colon\msc{T} \acts \cQ \otimes \msc{B}$ and of $\dot{\tau}\colon\msc{T} \acts \cQ \otimes \msc{B}$ into $\dot{\sigma}\colon\msc{S}\acts \cQ \otimes \msc{A}$.
    
    We say that $\ti{\sigma}\colon\msc{S}\acts \msc{A}$ and  $\ti{\tau}\colon\msc{T} \acts \msc{B}$ are \emph{strictly mutually embeddable} (and we write $\ti{\sigma}\colon\msc{S}\acts \msc{A} \sim_{s \cM \cE} \ti{\tau}\colon\msc{T} \acts \msc{B}$) if each can be strictly embedded into the other. If, in addition, the strict embeddings can be chosen to be of finite index, then we write $\ti{\sigma}\colon\msc{S}\acts \msc{A} \sim_{s \cM \cE_{FI}} \ti{\tau}\colon\msc{T} \acts \msc{B}$. Given an Abelian group $\cQ$, we write $\ti{\sigma}\colon\msc{S}\acts \msc{A} \sim_{s \cM \cE_{\cQ}} \ti{\tau}\colon\msc{T} \acts \msc{B}$ if $\dot{\sigma}\colon\msc{S}\acts \cQ \otimes \msc{A} \sim_{s \cM \cE} \dot{\tau}\colon\msc{T} \acts \cQ \otimes \msc{B}$. 
    
    We write $\ti{\sigma}\colon\msc{S}\acts \msc{A} \cong_{\cQ} \ti{\tau}\colon\msc{T} \acts \msc{B}$, and call $\ti{\sigma}\colon\msc{S}\acts \msc{A}$ and $\ti{\tau}\colon\msc{T} \acts \msc{B}$ \emph{isomorphic over $\cQ$}, if there exists a full and strict embedding of $\dot{\sigma}\colon\msc{S}\acts \cQ \otimes \msc{A}$ into $\dot{\tau}\colon\msc{T} \acts \cQ \otimes \msc{B}$. 
    \item We say that the algebraic actions $\sigma\colon S\acts A$ and $\tau\colon T \acts B$ are \emph{isomorphic} if there is a pair $(\mft,\mfb)$, where $\mft\colon S\to T$ is an isomorphism of semigroups and $\mfb\colon A\to B$ is a group isomorphism such that $\mfb(\sigma_s(x))=\tau_{\mft(s)}(\mfb(x))$ for all $s\in S$ and $x\in A$.
\end{enumerate} 
\edefin

\bremark
The dual notion of a strict (algebraic) embedding in our sense is an algebraic factor map in the sense of \cite[Definition~8.1]{SchPIMS}.
\setlength{\parindent}{0.5cm} \setlength{\parskip}{0cm}

If $(\mft,\mfb)$ is a finite index embedding of $\ti{\sigma}\colon\msc{S}\acts \msc{A}$ into $\ti{\tau}\colon\msc{T} \acts \msc{B}$ and $\msc{A}$ and $\msc{B}$ are torsion-free, then $\msc{A}$ and $\msc{B}$ are quasi-isomorphic in the sense of \cite[Definition~3.3]{Thomas}.

If $\ti{\sigma}\colon\msc{S}\acts \msc{A} \sim_{s \cM \cE_{FI}} \ti{\tau}\colon\msc{T} \acts \msc{B}$, then we also call $\ti{\sigma}\colon\msc{S}\acts \msc{A}$ and $\ti{\tau}\colon\msc{T} \acts \msc{B}$ \emph{finitely algebraically equivalent} (compare \cite[Definition~8.1]{SchPIMS}).

We will mostly be interested in our notions involving an Abelian group $\cQ$ when $\cQ = \Qz$.

If $\msc{A}$ and $\msc{B}$ are torsion-free, then $\sim_{\cM \cE_{FI}}$ implies $\sim_{\cM \cE_{\Qz^{\cong}}}$ and $\sim_{s \cM \cE_{FI}}$ implies $\cong_{\Qz}$.

If $\msc{A}$ and $\msc{B}$ are torsion-free and of finite rank, then $\sim_{\cM \cE}$ implies $\sim_{\cM \cE_{FI}}$ and $\sim_{s \cM \cE}$ implies $\sim_{s \cM \cE_{FI}}$ (see, e.g., \cite[Exercise~5 in \S~92]{Fuchs2}).
\eremark
\setlength{\parindent}{0cm} \setlength{\parskip}{0.5cm}

\bdefin
Given algebraic actions $\sigma \colon S \acts A$ and $\tau \colon T \acts B$, let \hypertarget{(I$_s$)}{(I$_s$)} be the symmetrized version of condition \hyperlink{(I)}{(I)} from Definition~\ref{def:I,II}, i.e., condition \hyperlink{(I)}{(I)} holds and the analogue of \hyperlink{(I)}{(I)} with reversed roles for $\sigma$  and $\tau$ holds as well. Similarly, let \hypertarget{(II$_s$)}{(II$_s$)}, \hypertarget{(III$_s$)}{(III$_s$)} and \hypertarget{($*_s$)}{($*_s$)} be the symmetrized versions of \hyperlink{(II)}{(II)}, \hyperlink{(III)}{(III)} and \hyperlink{(*)}{($*$)} from Definition~\ref{def:I,II}, Definition~\eqref{def:(III)}, and Definition~\eqref{def:STAR}.
\edefin

\bdefin
Suppose $\msc{A}$ is torsion-free and of finite rank. We say that $\ti{\sigma}\colon\msc{S}\acts\msc{A}$ is \emph{strongly faithful} if
\begin{equation}
\label{SF}\tag{SF}
    \dot{\sigma}_s=\rho\dot{\sigma}_t\rho^{-1} \text{ implies }s=t \text{ for all }s,t\in\msc{S}\text{ and } \rho\in\Aut(\Qz\otimes\msc{A}),
\end{equation}
where $\dot{\sigma}_s:=\id_{\Qz}\otimes\ti{\sigma}_s$.
\edefin
\bremark
\label{rem:SF}
If $\det\circ\dot{\sigma}\colon \msc{S} \to\Qz^\times$ is injective, then $\ti{\sigma}\colon\msc{S}\acts\msc{A}$ satisfies \eqref{SF}.
\eremark

We are now ready for the main result of this section.

\btheo
\label{thm:frmain}
Assume that $\sigma \colon S \acts A$ and $\tau \colon T \acts B$ are algebraic actions satisfying our standing assumptions from \S~\ref{ss:Standing}, with globalizations $\ti{\sigma}\colon\msc{S}\acts \msc{A}$ and $\ti{\tau}\colon\msc{T} \acts \msc{B}$. Suppose that $\msc{A}$ and $\msc{B}$ are torsion-free.
\setlength{\parindent}{0cm} \setlength{\parskip}{0cm}

\begin{enumerate}[\upshape(i)]
    \item If \hyperlink{(I$_s$)}{(I$_s$)} holds and there exists an isomorphism of topological groupoids
    \[
    (\msc{A} \rtimes \msc{S}) \ltimes \ol{A} \cong (\msc{B} \rtimes \msc{T}) \ltimes \ol{B},
    \]
    then $\ti{\sigma}\colon\msc{S}\acts \msc{A} \sim_{\cM \cE_{\Qz}} \ti{\tau}\colon\msc{T} \acts \msc{B}$.
    
    \item If \hyperlink{($*_s$)}{($*_s$)} holds and there exists an isomorphism of topological groupoids
    \[
    (S^{-1}A\rtimes S^{-1}S)\ltimes \ol{A}\cong (T^{-1}B\rtimes T^{-1} T)\ltimes\ol{B},
    \]
    then $\ti{\sigma}\colon S^{-1}S \acts S^{-1}A \sim_{\cM \cE} \ti{\tau}\colon T^{-1}T \acts T^{-1}B$.
\end{enumerate}
If, in addition, $\msc{A}$ and $\msc{B}$ have finite rank, then we obtain $\ti{\sigma}\colon\msc{S}\acts \msc{A} \sim_{\cM \cE_{\Qz^{\cong}}} \ti{\tau}\colon\msc{T} \acts \msc{B}$ in (i) and $\ti{\sigma}\colon S^{-1}S \acts S^{-1}A \sim_{\cM \cE_{FI}} \ti{\tau}\colon T^{-1}T \acts T^{-1}B$ in (ii).

If, in addition, $\msc{A}$ and $\msc{B}$ have finite rank and $\ti{\sigma}\colon\msc{S}\acts\msc{A}$, $\ti{\tau}\colon\msc{T}\acts\msc{B}$ both satisfy \eqref{SF}, then we obtain $\ti{\sigma}\colon\msc{S}\acts \msc{A} \cong_{\Qz} \ti{\tau}\colon\msc{T} \acts \msc{B}$ in (i) and $\ti{\sigma}\colon S^{-1}S \acts S^{-1}A \sim_{s \cM \cE_{FI}} \ti{\tau}\colon T^{-1}T \acts T^{-1}B$ (i.e., $\ti{\sigma}\colon S^{-1}S \acts S^{-1}A$ and $\ti{\tau}\colon T^{-1}T \acts T^{-1}B$ are finitely algebraically equivalent) in (ii).
\etheo
\setlength{\parindent}{0cm} \setlength{\parskip}{0cm}

\bproof
Everything except the last claim follows from Theorem~\ref{thm:mainthm}. For the last claim, suppose that $(\mft,\mfb)$ is an embedding of $\dot{\sigma}\colon\msc{S}\acts\Qz\otimes\msc{A}$ into $\dot{\tau}\colon\msc{T}\acts\Qz\otimes\msc{B}$ and $(\mfs,\mfa)$ is an embedding of $\dot{\tau}\colon\msc{T}\acts\Qz\otimes\msc{B}$ into $\dot{\sigma}\colon\msc{S}\acts\Qz\otimes\msc{A}$. For $s\in\msc{S}$, we have $(\mfa\circ\mfb) \dot{\sigma}_t(\mfa\circ\mfb)^{-1}=\mfa \dot{\tau}_{\mft(s)}\mfa^{-1}=\dot{\tau}_{\mfs\circ\mft(s)}$, so that \eqref{SF} implies $(\mfs\circ\mft)(s)=s$. Hence, $\mfs\circ\mft=\id_\msc{S}$, so by symmetry, $\mft\circ\mfs=\id_\msc{T}$. Our assumptions imply that $\Qz\otimes\msc{A}$ and $\Qz\otimes\msc{B}$ have the same dimension as rational vector spaces, so that the injective maps $\mfa$ and $\mfb$ are invertible. The second part of the last claim is similar using that any injective endomorphism of a torsion-free finite rank Abelian group necessarily has finite index image (see, e.g., \cite[Exercise~92.5]{Fuchs2}).
\eproof

\bremark
In some of our examples, $A$ will be finitely generated, in which case we take $\fA_k = A$ for all $k$ in condition \hyperlink{(*)}{($*$)}, so that the requirement $\sigma_s(\fA_k) \subseteq \fA_k$ in \hyperlink{(*)}{($*$)} is automatic.
\eremark
\setlength{\parindent}{0cm} \setlength{\parskip}{0.5cm}

\bremark
\label{rem:GPDCartanTFG}
Suppose that both $\sigma \colon S \acts A$ and $\tau \colon T \acts B$ are exact. Then the corresponding groupoids are effective and minimal by \cite[Theorem~4.14]{BL2} (see also \cite[Lemma~2.23]{KM}) and \cite[Corollary~7.4]{BL2}. Hence the following are equivalent:
\setlength{\parindent}{0cm} \setlength{\parskip}{0cm}

\begin{enumerate}[\upshape(i)]
    \item $(\msc{A} \rtimes \msc{S}) \ltimes \ol{A}$ and $(\msc{B} \rtimes \msc{T}) \ltimes \ol{B}$ are isomorphic as topological groupoids;
    
    \item $(\fA_{\sigma},\fD_{\sigma})$ and $(\fA_{\tau},\fD_{\tau})$ are isomorphic as Cartan pairs, where $\fA_{\sigma}$ and $\fA_{\tau}$ are as in \cite[Definition~3.1]{BL2}, and $\fD_{\sigma}$ and $\fD_{\tau}$ are as in \cite[Proposition~3.30]{BL2};
    
    \item $\bm{F}((\msc{A} \rtimes \msc{S}) \ltimes \ol{A})$ and $\bm{F}((\msc{B} \rtimes \msc{T}) \ltimes \ol{B})$ are isomorphic as abstract groups; 

    \item $\bm{D}((\msc{A} \rtimes \msc{S}) \ltimes \ol{A})$ and $\bm{D}((\msc{B} \rtimes \msc{T}) \ltimes \ol{B})$ are isomorphic as abstract groups.
\end{enumerate}
This follows from \cite{Ren} (see also \cite{Raad}) and \cite[Theorems~0.2 and 3.3]{Rub} (or \cite[Theorem~3.10]{Mat15}).
\eremark
\setlength{\parindent}{0cm} \setlength{\parskip}{0.5cm}

\section{Algebraic actions on finite rank torsion-free Abelian groups}

In this section, we apply our rigidity results to example classes of algebraic actions on finite rank torsion-free Abelian groups.

\subsection{Algebraic actions of torsion-free Abelian monoids whose dual actions are mixing}
\label{ss:mixing}

Let $\sigma \colon S \acts A$ be an algebraic action, with $A$ Abelian. Let $\hat{\sigma} \colon S \acts \widehat{A}$ be the dual action as in \cite[Remark~2.2]{BL2} and denote by $\mu$ the normalized Haar measure on $\widehat{A}$. Recall (see, for instance, \cite[\S~1]{Sch} or \cite[Definition~1.5]{Wal}) that $\hat{\sigma}$ is (strongly) mixing (with respect to $\mu$) if for all Borel subsets $X, Y$ of $\widehat{A}$ we have
\[
  \lim_{s \to \infty} \mu(X \cap \hat{\sigma}_s(Y)) = \mu(X) \mu(Y).
\]
If $S$ has no non-trivial finite subsemigroups, we have the following relation between the mixing property of $\hat{\sigma}$ and condition \eqref{F} for $\sigma$:
\bremark
\label{rem:Mixing=F}
Assume that $S$ has no non-trivial finite subsemigroups. Then $\hat{\sigma}$ is mixing if and only if we have, for all $0 \neq a \in A$ and $1 \neq s \in S$, that $\sigma_s(a) \neq a$, i.e., the analogue of condition \eqref{F} holds for $\sigma \colon S \acts A$.
\setlength{\parindent}{0.5cm} \setlength{\parskip}{0cm}

Indeed, in general (without our assumption on $S$), $\hat{\sigma}$ is mixing if and only if for all infinite subsemigroups $S' \subseteq S$ and $0 \neq a \in A$, we have that $\# \menge{\sigma_{s'}(a)}{s' \in S'} = \infty$ (see \cite[Theorem~1.6]{Sch} and also \cite{Ber}, in particular \cite[Theorem~2.1]{Ber}). It is now straightforward to see that, if $S$ has no non-trivial finite subsemigroups, the latter statement is equivalent to the condition that for all $0 \neq a \in A$, we have $\# \menge{s \in S}{\sigma_s(a) = a} < \infty$. This condition, in turn, is equivalent to the statement that we have $\sigma_s(a) \neq a$ for all $0 \neq a \in A$ and $1 \neq s \in S$ (again assuming that $S$ has no non-trivial finite subsemigroups), as desired, because $\menge{s \in S}{\sigma_s(a) = a}$  is always a subsemigroup of $S$.
\eremark

Note that the condition that $S$ has no non-trivial finite subsemigroups is in particular satisfied if $S$ is torsion-free, in the sense that for all $s_1, s_2 \in S$ and $i \in \Zz_{>0}$, $s_1^i = s_2^i$ implies that $s_1=s_2$.
\setlength{\parindent}{0cm} \setlength{\parskip}{0.5cm}

With the help of this observation, let us now present the first example class to which we can apply our general rigidity results.
\bcor
\label{cor:mixing}
Assume that $S$ and $T$ are non-trivial, Abelian, cancellative, torsion-free monoids, that $A$ and $B$ are torsion-free Abelian groups of finite rank, and that $\sigma \colon S \acts A$ and $\tau \colon T \acts B$ are non-automorphic faithful algebraic actions. Further suppose that there exist $s\in S$ and $t\in T$ such that the dual actions $\hat{\sigma}\vert_{s^\Nz}:s^\Nz\acts\widehat{A}$ and $\hat{\tau}\vert_{t^\Nz}:t^\Nz\acts\widehat{B}$ are mixing, where $s^\Nz:=\{s^n : n\in\Nz\}$ and $t^\Nz:=\{t^n : n\in\Nz\}$.
Let $\ti{\sigma} \colon S^{-1} S \acts S^{-1} A$ and $\ti{\tau} \colon T^{-1} T \acts T^{-1} B$ be the canonical globalizations as in \cite[Example~2.4]{BL2}. 

If there exists an isomorphism of topological groupoids
    \[
    (S^{-1}A\rtimes S^{-1}S)\ltimes \ol{A}\cong (T^{-1}B\rtimes T^{-1} T)\ltimes\ol{B},
    \]
    then 
$\ti{\sigma}\colon S^{-1}S \acts S^{-1}A \sim_{\cM \cE_{FI}} \ti{\tau}\colon T^{-1}T \acts T^{-1}B$.
\ecor
\setlength{\parindent}{0cm} \setlength{\parskip}{0cm}

\bproof
First of all, note that condition \eqref{JF} is satisfied because of \cite[Proposition~7.5]{BL2} and \eqref{FI} holds by \cite[Example~7.6]{BL2}. Thus $\sigma$ and $\tau$ satisfy our standing assumptions from \S~\ref{ss:Standing}. Moreover, it is straightforward to check that $S^{-1} S$ and $T^{-1} T$ are torsion-free and that $\rk_\Zz S^{-1} A = \rk_\Zz A$, $\rk_\Zz T^{-1} B = \rk_\Zz B$. Now our statement follows from Theorem~\ref{thm:frmain}~(ii) for the finite rank case because of Remark~\ref{rem:Mixing=F}.
\eproof
\setlength{\parindent}{0cm} \setlength{\parskip}{0.5cm}

Let us briefly explain the conclusion of our results for the case of (duals of) toral endomorphisms.
\bex
\label{ex:ToralEndo}
Let $a\in \M_n(\Zz)$ and $b\in \M_m(\Zz)$ with $|\det(a)|,|\det(b)|>1$, where $n,m\in\Zz_{>0}$. If $a$ and $b$ both have no roots of unity as eigenvalues, then the duals of the $\Nz$-actions $\sigma\colon\Nz\acts \Zz^n$ and $\tau\colon\Nz\acts \Zz^m$ given by $\sigma_k(v)=a^kv$ and $\tau_k(w)=b^kw$ for $k\in\Nz$, $v\in\Zz^n$, and $w\in\Zz^m$ are mixing. 
Suppose that the corresponding groupoids are isomorphic. Then, since \eqref{SF} holds in this case, Theorem~\ref{thm:frmain} implies that $n=m$ and that the matrices $a$ and $b$ must be conjugate over $\Qz$, i.e., there exists $c\in\GL_n(\Qz)$ such that $a=cbc^{-1}$.
\eex

\subsection{Canonical endomorphisms of torsion-free finite rank Abelian groups}
\label{ss:endoabgroups}

Let $A\subseteq \Qz^n$ be a torsion-free Abelian group of rank $n\in\Zz_{>0}$. The multiplicative monoid $\Zz^\times:=\Zz\setminus\{0\}$ acts on $A$ by multiplication: Each $s\in \Zz^\times$ gives rise to the endomorphism $\sigma_s\colon A\to A$ given by $\sigma_s(x)= sx$. For any submonoid $M\subseteq\Zz^\times$, the associated algebraic action $\sigma^M\colon M\acts A$, where $\sigma^M:=\sigma\vert_M$, is faithful and satisfies \eqref{FI}. It is easy to see that $\sigma^M\colon M\acts A$ is non-automorphic if and only if there exists $m\in M$ such that $A$ is not $m$-divisible, i.e., $mA\subsetneq A$. Since $\det(\dot{\sigma}_s)=s^n$, we see that $\det\circ \dot{\sigma} \colon\Zz_{>0}\to\Qz^\times$ is injective, so that $\sigma^M\colon M\acts A$ satisfies \eqref{SF} (see Remark~\ref{rem:SF}).
Since $M$ is Abelian, we obtain a globalization $\ti{\sigma}^M\colon \gp{M}\acts M^{-1}A$, where $\gp{M}:=M^{-1}M\subseteq\Qz^\times$ acts on $M^{-1}A:=\bigcup_{s\in M}\frac{1}{s}A$ by multiplication. It is easy to see that $\ti{\sigma}^M\colon \gp{M}\acts M^{-1}A$ satisfies \eqref{F}. For $s,s'\in\Zz_{>0}$, we have $sA\cap s'A=\lcm(s,s')A$ (see, e.g., \cite[\S~20]{Fuchs1}). Hence, the $M$-constructible subgroups for $M\acts A$ are given by $\{sA : s\in M \}$. From this, we see that $\ti{\sigma}^M\colon \gp{M}\acts M^{-1}A$ satisfies \eqref{JF}.
Let $x_1,...,x_n\in A$ be rationally independent, and put $\fA:=\spn_\Zz(\{x_1,...,x_n\})\cong\Zz^n\subseteq A$. For each $k\in\Zz_{>0}$, let
\[
\fA_k:=\{x\in A : (k!)x\in\fA\}=A\cap (k!)^{-1}\fA.
\]
Each $\fA_k$ is an $\Zz^\times$-invariant finitely generated subgroup of $A$, and we have $A=\bigcup_k\fA_k$.

We can now apply Theorem~\ref{thm:frmain} to obtain the following result:

\bcor
Let $A$ and $B$ be torsion-free finite rank Abelian groups. Let $M,N\subseteq\Zz^\times$ be submonoids such that there exist $m\in M$ and $n\in N$ with $mA\subsetneq A$ and $nB\subsetneq B$. If there is an isomorphism of topological groupoids 
\[
(M^{-1}A\rtimes \gp{M})\ltimes\ol{A}\cong (N^{-1}B\rtimes \gp{N})\ltimes\ol{B},
\]
then the actions $\gp{M}\acts M^{-1}A$ and $\gp{N}\acts N^{-1}B$ are finitely algebraically equivalent.
\ecor
\bremark
If $\Zz_{>0}\subseteq M$, then $M\acts A$ is exact if and only if the Ulm  subgroup of $A$ vanishes (cf. \cite[\S~1.6]{Fuchs1}).
\eremark

\subsection{Actions adding scalars to algebraic actions of groups}
\label{ss:addingscalars}
Let $\Gamma\subseteq\SL_n(\Zz)$ be any subgroup, and let $M\subseteq \Zz_{>0}$ a submonoid. Then, $M\Gamma:=\{a\gamma :a\in M,\gamma\in\Gamma\}$ is a submonoid of $\M_n(\Zz)^\times:=\{x\in\M_n(\Zz) : \det(x)\neq 0\}$, where we view $M$ as a submonoid of $\M_n(\Zz)$ via the diagonal embedding. Since $\M_n(\Zz)^\times$ acts canonically on $\Zz^n$, we obtain a faithful algebraic action $M\Gamma\acts\Zz^n$. It is easy to see that $M\Gamma\acts\Zz^n$ is exact if and only if $M$ is non-trivial.
Note that $\gp{M}\Gamma\acts (M^{-1}\Zz)^n$ is a globalization for $M\Gamma\acts\Zz^n$ that satisfies \eqref{JF}.  Since $\Gamma$ acts by automorphisms on $\Zz^n$ that commute with the action of $M$, we have $\cC_{M\Gamma\acts \Zz^n}=\cC_{M\acts \Zz^n}$. Let $\ol{\Zz}_M^n$ denote the completion of $\Zz^n$ with respect to the family $\cC_{M\acts \Zz^n}$.

\bremark
The globalization $\gp{M}\Gamma\acts (M^{-1}\Zz)^n$ often will not satisfy \eqref{F}. For instance, take $M=
\Zz_{>0}$ and $\Gamma=\SL_2(\Zz)$. Then $\id-\gamma$ is not injective on $M^{-1}\Zz^2=\Qz^2$, where $\gamma=\bigl(\begin{smallmatrix}1 & 1 \\0 & 1\end{smallmatrix}\bigl)\in\SL_2(\Zz)$.
\eremark

In order to apply our rigidity result in this setting, we need an observation on subgroups of $\SL_n(\Zz)$, which comes from \cite[Example~26.8~\&~Lemma~26.16]{DK}.
\blemma
\label{lem:Gammaprops}
Let $\Gamma\subseteq\SL_n(\Zz)$ be a subgroup. If $\gamma\alpha\gamma^{-1}=\alpha^{\kappa}$ for $\alpha,\gamma\in\Gamma$ and $\kappa\in\Zz_{>1}$, then $\alpha\in\tor(\Gamma)$.
\elemma
\setlength{\parindent}{0cm} \setlength{\parskip}{0cm}

\bproof
Suppose $\gamma\alpha\gamma^{-1}=\alpha^{\kappa}$ for $\alpha,\gamma\in\Gamma$ and $\kappa\in\Zz_{>1}$. Let $p$ be a prime divisor of $\kappa$. For $l\geq 1$, consider the congruence subgroup $\Gamma(p^l):=\{a\in\SL_n(\Zz) : a \equiv I_n \mod p^l\}$, and put $\Gamma_{p^l}:=\Gamma\cap\Gamma(p^l)$. Since $\Gamma_p$ is a finite index subgroup of $\Gamma$, we can find $m\in\Zz_{>0}$ such that $\alpha^m\in\Gamma_p$. If $\alpha^m\neq I_n$, then since $\bigcap_l\Gamma_{p^l}=\{I_n\}$, we can find $l\geq 1$ with $\alpha^m\notin\Gamma_{p^l}$. Now $\alpha^m\Gamma_{p^l}$ and $\alpha^{\kappa m}\Gamma_{p^l}$ have the same order in the $p$-group $\Gamma_p/\Gamma_{p^l}$ because $\gamma\alpha\gamma^{-1}=\alpha^{\kappa}$, which is a contradiction since $p\mid \kappa$.
\eproof
\setlength{\parindent}{0cm} \setlength{\parskip}{0.5cm}

\bcor
\label{cor:GammaLambda}
Let $\Gamma,\Lambda\subseteq \SL_n(\Zz)$ be subgroups and $M,N\subseteq\Zz_{>0}$ nontrivial submonoids such that for every $\gamma\in\tor(\Gamma)$, there exists $s\in M$ with $\gcd(\ord(\gamma),s)>1$, and for every $\lambda\in\tor(\Lambda)$, there exists $t\in N$ with $\gcd(\ord(\lambda),t)>1$. If there is an isomorphism of topological groupoids
\[ 
(M^{-1}\Zz^n\rtimes \gp{M}\Gamma)\ltimes\ol{\Zz}^n_M\cong (N^{-1}\Zz^m\rtimes\gp{N}\Lambda)\ltimes\ol{\Zz}^m_N,
\]
then $M=N$ and there exist $g,h\in\M_n(N^{-1}\Zz)\cap\GL_n(\Qz)$ such that $\Gamma\subseteq g\Lambda g^{-1}$ and $h\Lambda h^{-1}\subseteq \Gamma$.
\ecor
\setlength{\parindent}{0cm} \setlength{\parskip}{0cm}

\bproof
By Lemma~\ref{lem:Gammaprops}, condition \hyperlink{($*_s$)}{($*_s$)} holds, so this follows from Theorem~\ref{thm:frmain}.
\eproof

\bremark
The assumptions on $\Gamma,\Lambda$ and $M,N$ in the statement Corollary~\ref{cor:GammaLambda} are satisfied, for instance, if $M=N=\Zz_{>0}$ or if $\Gamma$ and $\Lambda$ are torsion-free (and $M,N\neq \{1\}$). 
\eremark

\bremark
If in the statement of Corollary~\ref{cor:GammaLambda}, $\Gamma$ and $\Lambda$ are not conjugate via an element in $\GL_n(\Qz)$ to any of their proper subgroups, then the conclusion can be strengthened to the following: $M=N$ and there exists $g\in\M_n(N^{-1}\Zz)\cap\GL_n(\Qz)$ such that $\Gamma=g\Lambda g^{-1}$. This holds, for instance, if $\Gamma$ and $\Lambda$ are co-Hopfian.
\eremark
\setlength{\parindent}{0cm} \setlength{\parskip}{0.5cm}

\subsection{Arithmetic dynamical systems}
\label{ss:arithmeticdynamics}
Let us consider the algebraic $\Nz$-actions studied by Chothi, Everest, and Ward in \cite{CEW}.
Let $K$ be a number field with ring of integers $R$, and let $\cP_K$ denote the set of non-zero prime ideals of $R$. For $\p\in \cP_K$, let $v_\p$ and $|\cdot|_\p$ denote the associated additive and multiplicative $\p$-adic valuations on $K$, respectively. Given a subset $\cS\subseteq \cP_K$, the corresponding ring of \emph{$\cS$-integers} is
\[
R_{\cS}:=\{x\in K : |x|_\p\leq 1\text{ for every }\p\in \cP_K\setminus\cS\},
\]
i.e., $R_{\cS}$ consists of the elements of $K$ that are $\p$-adic integers for every $\p\in \cP_K\setminus\cS$. For $\xi\in R_{\cS}^\times=R_{\cS}\setminus\{0\}$, the map $m_\xi\colon R_{\cS}\to R_{\cS}$ given by $m_\xi(x)=\xi x$ is an injective endomorphism of the additive group of $R_{\cS}$. The dual action $\widehat{m}_\xi\colon\Nz\acts\widehat{R_{\cS}}$ is called an \emph{arithmetic $\cS$-integer dynamical system}, see \cite[\S~2]{CEW}. The group of units (i.e., invertible elements) in $R_{\cS}$ is $R_{\cS}^*=\{x\in K^* : |x|_\p=1\text{ for every }\p\in \cP_K\setminus\cS\}$. Note that $R_{\cS}$ is a proper subring of $K$ if and only if $\cS\subsetneq \cP_K$. Also note that $R_{\cS}^*\subsetneq R_{\cS}^\times$ whenever $R_{\cS}\subsetneq K$. Let us record some basic observations about the algebraic action $m_\xi\colon\Nz\acts R_{\cS}$. Let $\ti{m}_\xi\colon R_{\cS}[1/\xi]\to R_{\cS}[1/\xi]$ be given by $\ti{m}_\xi(x)=\xi x$. Then $\ti{m}_\xi\colon\Zz\acts R_{\cS}[1/\xi]$ is a globalization of $m_\xi\colon\Nz\acts R_{\cS}$.
\blemma
\begin{enumerate}[\upshape(i)]
    \item $m_\xi\colon\Nz\acts R_{\cS}$ is faithful if and only if $\xi$ is not a root of unity.
    \item $m_\xi\colon\Nz\acts R_{\cS}$ is exact if and only if it is non-automorphic if and only if $\xi$ is a non-unit.
    \item $\ti{m}_\xi\colon\Zz\acts R_{\cS}[1/\xi]$ satisfies \eqref{JF}.
\end{enumerate}
\elemma
\setlength{\parindent}{0cm} \setlength{\parskip}{0cm}

\bproof
(i) and (iii) are obvious. For (ii), let $\xi\in R_{\cS}^\times\setminus R_{\cS}^*$. Then there exists $\p\in\cP_K\setminus\cS$ such that $|\xi|_\p<1$. If $x\in R_{\cS}$ lies in $\bigcap_{n\geq 0}\xi^nR_{\cS}$, then for each $n\geq 0$, we can write $x=\xi^n y_n$ for some $y_n\in R_{\cS}$.  Now we have $|x|_\p=|\xi|_\p^n|y_n|_\p\leq |\xi|_\p^n$ for every $n$, so that $x=0$. The other implications in (ii) are easy to see.
\eproof
\setlength{\parindent}{0cm} \setlength{\parskip}{0.5cm}

\bremark
An element $x\in K$ is integral over $\Zz$ if and only if $\Zz[x]$ is finitely generated as a $\Zz$-module, so the additive group of $R_{\cS}$ is not finitely generated whenever $\cS\neq\emptyset$.
\eremark

For each $k\in\Zz_{>0}$, let $\fA_k:=R_{\cS} \cap (k!)^{-1}R=\{x\in R_{\cS} : (k!)x\in R\}$.
Then $R_{\cS}=\bigcup_k \fA_k$, and every $\fA_k$ is finitely generated and invariant under $R^\times$.

\bcor
\label{cor:arithmeticdynamics}
Let $K_1$ and $K_2$ be number fields with rings of algebraic integers $R_1$ and $R_2$, respectively, let $\cS\subsetneq\cP_{K_1}$ and $\cT\subsetneq \cP_{K_2}$ by proper subsets of primes, and let $\xi\in R_1^\times\setminus R_{1,\cS}^*$ and $\eta\in R_2^\times\setminus R_{2,\cT}^*$.
If there is an isomorphism of topological groupoids 
\[
(R_{1,\cS}[1/\xi]\rtimes \gp{\xi})\ltimes\ol{R_{1,\cS}}\cong (R_{2,\cT}[1/\eta]\rtimes \gp{\eta})\ltimes\ol{R_{2,\cT}},
\]
then $\tilde{m}_{\xi}\acts R_{1,\cS}[1/\xi]$ and $\tilde{m}_{\eta}\acts R_{2,\cT}[1/\eta]$ are finitely algebraically equivalent. 
\ecor
\setlength{\parindent}{0cm} \setlength{\parskip}{0cm}

\bproof
Condition \hyperlink{($*_s$)}{($*_s$)} is satisfied, so Theorem~\ref{thm:frmain} yields the result.
\eproof

\bremark
Given any $\xi\in R_{\cS}^\times\setminus R_{\cS}^*$, there exists $l\in\Zz_{>0}$ such that $l\xi\in R$, and then the pair $(\id,m_l)$ is a strict, finite index embedding of $m_\xi\colon \Nz\acts R_{\cS}$ into $m_{l\xi}\colon \Nz\acts R_{\cS}$; in particular, $m_\xi\colon \Nz\acts R_{\cS}$ and $m_{l\xi}\colon \Nz\acts R_{\cS}$ are isomorphic over $\Zz[l^{-1}]$. Thus, up to inverting an integer, Corollary~\ref{cor:arithmeticdynamics} applies to all (faithful, exact) actions of the form $m_\xi\colon \Nz\acts R_{\cS}$.
\eremark
\setlength{\parindent}{0cm} \setlength{\parskip}{0.5cm}

\section{Algebraic actions from rings}
\label{s:NT}

The \emph{rank} of a ring $R$ is defined to be the rank of the additive group of $R$, i.e., the dimension of $\Qz\otimes_\Zz R$ as a vector space over $\Qz$. We shall say that $R$ is \emph{torsion-free} if the additive group of $R$ is a torsion-free (Abelian) group.
Examples of torsion-free rings of finite rank include integral group rings of finite groups and $R^n$ or $\M_n(R)$, where $R$ is an order in a central simple algebra over an algebraic number field. 

\subsection{General preparations}
Let $R$ be a unital torsion-free ring of finite rank $n\in\Zz_{>0}$. Then $\Qz\otimes_\Zz R$ is an $n$-dimensional $\Qz$-algebra containing $R$ as a full subring. The sum of two elementary tensors in $\Qz\otimes_\Zz R$ is again an elementary tensor, so that $\Qz\otimes_\Zz R=\Qz R:=\{q\otimes x : q\in \Qz,x\in R\}$. Moreover, if $\cL\subseteq R$ is a full rank subgroup, then $\Qz\otimes_\Zz\cL=\Qz\cL$, and $\Qz \cL = \Qz R$. 
Each $a\in \Qz R$ gives rise to a $\Qz$-linear map $\dot{\sigma}_a\colon \Qz R\to\Qz R$ given by $x\mapsto ax$. We let  $\chi_a(t)$ denote the characteristic polynomial of this map, and put $N(a):=|\det(\dot{\sigma}_a)|$. For $a\in R^\times$, put $\sigma_a:=\dot{\sigma}_a\vert_{R}$.

Following the notation from \cite{Li:Ring}, we let $ R^\times$ denote the multiplicative monoid of left regular elements in $ R$, i.e., $R^\times$ consists of those $a\in R$ such that $\sigma_a$ is injective. Since $R$ is a torsion-free ring of finite rank, the element $a\in R$ is left regular if and only if $N(a)\neq 0$. The action of any submonoid $M\subseteq R^\times$ on (the additive group of) $R$ by left multiplication is faithful, by injective endomorphisms, and satisfies \eqref{FI} (see, e.g., \cite[Exercise~92.5]{Fuchs2}). If $\cL$ is a full rank additive subgroup of $R$ that is invariant under the action of $M$, then $M$ also acts faithfully on $\cL$ by injective endomorphisms and the action $M\acts\cL$ satisfies \eqref{FI}. Under the canonical inclusion $R\subseteq \Qz R$, $R^\times$ is carried into $(\Qz R)^*$. If $M\subseteq R^\times$ is a submonoid, then we let $\langle M\rangle$ denote the subgroup of $(\Qz R)^*$ generated by $M$. Since $\cL$ is of full rank, we have $\cL\cap R^\times\neq \emptyset$ and thus $M\acts \cL$ is faithful.

\bprop
\label{prop:extension}
For $i=1,2$, let $R_i$ be a torsion-free ring of rank $n$ and $M_i\subseteq R_i^\times$ a submonoid. 
Let $\cL$ be a rank $n$ subgroup of $R_1$, and assume that $\spn_\Zz(M_1)$ has finite index in $R_1$. 
If there is an injective additive group homomorphism $\mfb \colon \Qz\cL\to \Qz R_2$ and a group homomorphism $\mft \colon \gp{M_1}\to \langle M_2\rangle$ such that $\mfb(ax)=\mft(a)\mfb(x)$ for all $a\in M_1$ and $x\in \Qz\cL$, then there exists a unital $\Qz$-algebra isomorphism $\varphi\colon \Qz R_1\to \Qz R_2$ such that $\varphi\vert_{\langle M_1\rangle }=\mft$.
\eprop
\setlength{\parindent}{0cm} \setlength{\parskip}{0cm}

\bproof
First, we show that $\mfb(1)$ is invertible in $\Qz R_2$. For every $a\in M_1$, we have $\mfb(a)=\mfb(a1) = \mft(a)\mfb(1)$, so that $\mfb(x)$ lies in the $\Qz$-vector space $(\Qz R_2)\mfb(1)$ for all $x\in \spn_\Zz(M_1)$. 
Since $\mfb$ is injective and $\rk_\Zz(\spn_\Zz(M_1))=n$, we have $n=\rk_\Zz(\im(\mfb))\leq\dim_\Qz( (\Qz R_2)\mfb(1) ) \leq \dim_\Qz \Qz R_2=n$. Hence, $\dim_\Qz( (\Qz R_2)\mfb(1) ) = \dim_\Qz (\Qz R_2)$, which implies that $\mfb(1)$ is invertible.

We now define $\varphi\colon \Qz  R_1\to \Qz R_2$ by $\varphi(x):=\mfb(x)\mfb(1)^{-1}$. Clearly, $\varphi$ is additive and $\varphi(1)=1$. Since $\mfb$ is injective and $\mfb^{-1}$ is invertible, we see that $\varphi$ is also injective. Let $a\in M_1$. We obtain $\mft(a)=\mfb(a)\mfb(1)^{-1}=\varphi(a)$. Thus, for $a\in M_1$ and $x\in  R_1$, we have
\[
\varphi(ax)=\mfb(ax)\mfb(1)^{-1} = \mft(a)\mfb(x)\mfb(1)^{-1}=\mft(a)\varphi(x)=\varphi(a)\varphi(x).
\]
Set $\Mult(\varphi):=\{a\in  R_1 : \varphi(ax) = \varphi(a) \varphi(x) \text{ for all } x\in  R_1\}$. It is straightforward to see that $\Mult(\varphi)$ is a subring of $R_1$ containing $\Zz$ and $M_1$. Since $\spn_\Zz(M_1)$ is of finite index in $R_1$, it follows that $\Mult(\varphi)= R_1$. Hence $\varphi$ is a ring homomorphism. Thus $\varphi$ is a $\Qz$-algebra isomorphism $\Qz R_1\to\Qz R_2$ satisfying $\varphi\vert_{\gp{M_1}}=\mft$.
\eproof
\setlength{\parindent}{0cm} \setlength{\parskip}{0.5cm}

Given a torsion-free ring $R$ of finite rank, we let $\cO$ denote the integral closure of $\Zz$ in $\Qz R$. If $R$ is finitely generated, then $R\subseteq\cO$ by \cite[Theorem~1.10]{Reiner}. However, $\cO$ may not be a subring if $R$ is non-commutative.
Given a submonoid $M\subseteq  R^\times$, let $\widetilde{M}:=\langle M\rangle \cap \cO$. Note that $\widetilde{M}$ need not be closed under multiplication.

The following Corollary demonstrates criteria under which we can deduce rigidity results.

\bcor
\label{cor:QAlgIsom}
For $i=1,2$, suppose $R_i$ is a finitely generated torsion-free ring of finite rank and that $M_i\subseteq  R_i^\times$ is submonoid. Assume there exist $\Qz$-algebra isomorphisms $\varphi_1\colon \Qz R_1\to \Qz R_2$ and $\varphi_2\colon\Qz R_2\to\Qz R_1$ such that $\varphi_1(\langle M_1\rangle)\subseteq \langle M_2\rangle$ and $\varphi_2(\langle M_2\rangle)\subseteq \langle M_1\rangle$. If
\setlength{\parindent}{0cm} \setlength{\parskip}{0cm}

\begin{enumerate}
	\item[\textup{(S')}] $\psi(\langle M_1\rangle )\subseteq \langle M_1\rangle \implies \psi(\langle M_1\rangle )=\langle M_1\rangle$ for all $\psi\in \Aut_{\Qz\textup{-alg}}(\Qz R_1)$,
\end{enumerate}
then $\varphi_1(\langle M_1\rangle)=\langle M_2\rangle$ and $\varphi_2(\langle M_2\rangle)=\langle M_1\rangle$, so that $\langle M_1\rangle\acts \Qz R_1$ and $\langle M_2\rangle\acts \Qz R_2$ are isomorphic. 
\setlength{\parindent}{0cm} \setlength{\parskip}{0.5cm}

If
\setlength{\parindent}{0cm} \setlength{\parskip}{0cm}

\begin{enumerate}
	\item[\textup{(N)}] $M_i=\widetilde{M_i}$ (for $i=1,2$), and
	\item[\textup{(S)}] $\psi(M_1)\subseteq M_1\implies \psi(M_1)=M_1$ for all $\psi\in \Aut_{\Qz\textup{-alg}}(\Qz R_1)$,
\end{enumerate}
then $\varphi_1(M_1)=M_2$ and $\varphi_2(M_2)=M_1$; therefore, $M_1\acts \spn_\Zz(M_1)$ and $M_2\acts \spn_\Zz(M_2)$ are isomorphic, and, if each $\cO_i$, the integral closure of $\Zz$ in $\Qz R_i$, is closed under addition and $M_i$-invariant, then $M_1\acts\cO_1$ and $M_2\acts\cO_2$ are isomorphic.
\ecor
\setlength{\parindent}{0cm} \setlength{\parskip}{0cm}

\bproof
Let $\psi:=\varphi_2\circ\varphi_1\in \Aut_{\Qz\textup{-alg}}(\Qz R_1)$. Then $\psi(\gp {M_1})=\varphi_2(\varphi_1(\gp{M_1}))\subseteq \varphi_2(\gp{M_2})\subseteq \gp{M_1}$, so that $\psi(\gp{M_1})=\gp{M_1}$ by (S'). Now we have $\gp{M_1}=\psi(\gp{M_1})=\varphi_2(\varphi_1(\gp{M_1}))\subseteq \varphi_2(\gp{M_2})\subseteq \gp{M_1}$, so that $\varphi_2(\gp{M_2})=\gp{M_1}$.
\setlength{\parindent}{0.5cm} \setlength{\parskip}{0cm}

Now assume that (N) and (S) hold. Since $\Qz$-algebra homomorphisms preserve integrality, we have $\varphi_1(\cO_1)=\cO_2$ and $\varphi_2(\cO_2)=\cO_1$; moreover, we have $M_i\subseteq \cO_i$ for $i=1,2$, so our assumption that $\varphi_1(M_1)\subseteq\langle M_2\rangle$ and $\varphi_2(M_2)\subseteq\langle M_1\rangle$ forces $\varphi_1(M_1)\subseteq \cO_2\cap\langle M_2\rangle$ and $\varphi_2(M_2)\subseteq \cO_1\cap\langle M_1\rangle$.
We have 
\begin{equation*}
\psi(M_1)=\varphi_2(\varphi_1(M_1))\subseteq\varphi_2(\cO_2\cap\langle M_2\rangle)\overset{(N)}{=}\varphi_2(M_2) \subseteq \cO_2\cap \langle M_1\rangle\overset{(N)}{=}M_1,
\end{equation*}
so that condition (S) forces $\psi(M_1)=M_1$, so that $\varphi_1(M_1)=M_2$ and $\varphi_2(M_2)=M_1$.
\eproof
\setlength{\parindent}{0cm} \setlength{\parskip}{0.5cm}

\subsection{Groupoid rigidity when the acting monoid is Abelian}	
In this section, we specialise to the case where the acting monoids are Abelian. The following is an immediate consequence of Theorem~\ref{thm:mainthm} and Corollary~\ref{cor:QAlgIsom}.
\btheo
\label{thm:MactsL}
For $i=1,2$, suppose $R_i$ is a torsion-free finitely generated ring, $M_i\subseteq  R_i^\times$ an Abelian submonoid such that $\spn_\Zz(M_i)$ has finite index in $R_i$, and $\cL_i\subseteq  R_i$ an $M_i$-invariant full rank subgroup. Assume that there exists $a\in M_1$ such that $\cL_1\to\cL_1, x\mapsto (1-a)x$ is injective, and similarly for $M_2$.

If there is an isomorphism of topological groupoids $(M_1^{-1}\cL_1\rtimes \langle M_1\rangle)\ltimes \ol{\cL}_1\cong (M_2^{-1}\cL_2\rtimes \langle M_2\rangle)\ltimes \ol{\cL}_2$, then there exist $\Qz$-algebra isomorphisms $\varphi_1\colon \Qz R_1\overset{\cong}{\to}\Qz R_2$ and $\varphi_2\colon \Qz R_2\overset{\cong}{\to}\Qz R_1$ such that $\varphi_1(M_1)\subseteq \widetilde{M_2}$ and $\varphi_2(M_2)\subseteq \widetilde{M_1}$.
\etheo

We obtain the following rigidity results. 
\bcor
\label{cor:Mispan(Mi),Oi}
Suppose that, in addition to the assumptions in Theorem~\ref{thm:MactsL}, conditions (N) and (S) from Corollary~\ref{cor:QAlgIsom} hold and that $\cL_i=\spn_\Zz(M_i)=R_i$ or $\cL_i=\cO_i$, then the following statements are equivalent:
\setlength{\parindent}{0cm} \setlength{\parskip}{0cm}

\begin{enumerate}[\upshape(i)]
	\item the algebraic actions $M_1\acts R_1$ and $M_2\acts R_2$ are isomorphic;
	\item $(M_1^{-1}R_1\rtimes \langle M_1\rangle)\ltimes \ol{\cL}_1$ and  $(M_2^{-1} R_2\rtimes \langle M_2\rangle)\ltimes \ol{\cL}_2$ are isomorphic as topological groupoids.
\end{enumerate}
\ecor
\setlength{\parindent}{0cm} \setlength{\parskip}{0cm}

Note that $\cL_i=\cO_i$ requires that $\cO_i$ is closed under addition and $M_i$-invariant, neither of which is automatic.
\setlength{\parindent}{0cm} \setlength{\parskip}{0.5cm}

\subsubsection{Connection to Bhargava's work}

Torsion-free commutative rings with finitely generated additive groups have received a great deal of attention recently \cite{Bh1,Bh2,Bh3,Bhd1,Bh4,Bhd2}. 

\btheo
\label{thm:Bhargava}
Let $R_i$, $i=1,2$, be finitely generated torsion-free rings. If the groupoids $(\Qz R_1\rtimes (\Qz R_1)^*)\ltimes\ol{R}_1$ and $(\Qz R_2\rtimes (\Qz R_2)^*)\ltimes\ol{R}_2$ are isomorphic, then $\Qz R_1\cong \Qz R_2$ as $\Qz$-algebras.
\etheo
\setlength{\parindent}{0cm} \setlength{\parskip}{0cm}

\bproof
This follows from Theorem~\ref{thm:MactsL}, applied to $M_i = R_i\reg$ and $\cL_i = R_i$. Since $\Zz \subseteq R_i$, we just need to show that $\spn_\Zz(M_i) = R_i$. Indeed, for every $a \in R_i$, there exists $\kappa \in \Zz\reg$ such that $a + \kappa \in R_i\reg$: As $\sp(\dot{\sigma}_a)$ is finite, we have $0 \notin \sp(\dot{\sigma}_{a + \kappa}) = \sp(\dot{\sigma}_a + \kappa \, \id) = \sp(\dot{\sigma}_a) + \kappa$ for sufficiently big $\kappa$.
\eproof
\setlength{\parindent}{0cm} \setlength{\parskip}{0.5cm}

\subsubsection{Actions of congruence monoids on rings of algebraic integers}
\label{ssecNT}

Let $K$ be a number field with ring of integers $R$, and let $R_{\m,\Gamma}\subseteq R^\times=R\setminus\{0\}$ be a congruence monoid as in \cite[\S~3]{Bru1}, where $\m=\m_\infty\m_0$ is a modulus for $K$ and $\Gamma$ is a group of residues modulo $\m$. Let $C_\lambda^*(R\rtimes R_{\m,\Gamma})$ denote the left regular C*-algebra of the monoid $R\rtimes R_{\m,\Gamma}$ and $D_\lambda(R\rtimes R_{\m,\Gamma})$ the canonical Cartan subalgebra of $C_\lambda^*(R\rtimes R_{\m,\Gamma})$ (see \cite[\S~2.2]{BL}). Using the results from \cite[\S~2]{Bru1}, it is not difficult to show that the family of constructible subgroups for the multiplication action $R_{\m,\Gamma}\acts R$ is given by $\cC_{R_{\m,\Gamma}\acts R}=\{(0)\neq I\unlhd R : I \text{ is coprime with }\m\}$. In particular, the completion $\ol{R}$ of $R$ with respect to $\cC_{R_{\m,\Gamma}\acts R}$ depends only on the prime divisors of $\m$.

\blemma
\label{lem:congrumonoid1}
The ring $\spn_\Zz(R_{\m,\Gamma})$ is an order in $R$, i.e., it is of finite index in $R$.
\elemma
\setlength{\parindent}{0cm} \setlength{\parskip}{0cm}

\bproof
Since $R_{\m,\Gamma}$ contains $(1+\m_0)_+$, it follows that the subring of $R$ generated by $R_{\m,\Gamma}$ contains $(\m_0)_+$, the set of totally positive elements in the ideal $\m_0$. If $x\in\m_0$, choose $k\in\Nz\reg\cap\m_0$ such that $x+k$ is totally positive. Since $x=(x+k)-k$, we see that every element of $\m_0$ is a difference of totally positive elements. Hence, the ring generated by $R_{\m,\Gamma}$ contains $\m_0$, which implies that it is of finite index in $R$.
\eproof
\setlength{\parindent}{0cm} \setlength{\parskip}{0.5cm}

\blemma
\label{lem:congrumonoid2}
\begin{enumerate}[\upshape(i)]
    \item The monoid $R_{\m,\Gamma}$ satisfies conditions (N) and (S) as in Corollary~\ref{cor:QAlgIsom};
    \item the action $R_{\m,\Gamma}\acts R$ is exact.
\end{enumerate}
\elemma
\setlength{\parindent}{0cm} \setlength{\parskip}{0cm}

\bproof
We shall use the notation from \cite{Bru1}. (i):  First, let us show condition (N). By Proposition~\cite[Proposition~3.2]{Bru1}, we have
\[
\langle R_{\m,\Gamma}\rangle=\{x\in K^\times : v_\p(x)=0\text{ for all }\p\mid\m_0, [x]_\m\in\Gamma\},
\]
from which we see that $\langle R_{\m,\Gamma}\rangle\cap R=R_{\m,\Gamma}$. 

Now we verify that condition (S) holds. Let $\psi \in\Gal(K/\Qz)$. We have $\psi(R_{\m,\Gamma})=R_{\psi(\m),\psi(\Gamma)}$ where $\psi(\m)$ is the modulus defined by $w\mid\psi(\m)_\infty$ if and only if $w\circ\psi\mid \m_\infty$ and $\psi(\m)_0:=\psi(\m_0)$, and $\psi(\Gamma)$ is the image of $\Gamma$ under the isomorphism $(R/\m)^*\cong (R/\psi(\m))^*$ given by $[a]_\m\mapsto [\psi(a)]_{\psi(\m)}$.
Since $R_{\psi(\m),\psi(\Gamma)} = \psi(R_{\m,\Gamma}) \subseteq R_{\m,\Gamma}$, \cite[Proposition~9.2(1)]{Bru1} implies that $\psi(\m)\mid \m$, i.e., $\psi(\m)_0\mid\m_0$ and $\psi(\m)_\infty(w)=1\implies \m_\infty(w)=1$, and $\pi_{\m,\psi(\m)}(\psi(\Gamma))\subseteq \Gamma$, where $\pi_{\m,\psi(\m)}:(R/\m)^*\to (R/\psi(\m))^*$ is the canonical quotient map arising from the divisibility condition $\psi(\m)\mid \m$. Since $\psi(\m_0)$ and $\m_0$ have the same norm, $\psi(\m_0)\mid \m_0$ forces $\psi(\m_0)=\m_0$; since the finite sets $\supp(\psi(\m)_\infty)$ and $\supp(\m_\infty)$ have the same cardinality, $\supp(\psi(\m)_\infty)\subseteq \supp(\m_\infty)$ forces $\supp(\psi(\m)_\infty)=\supp(\m_\infty)$, i.e., $\psi(\m)_\infty=\m_\infty$.
Therefore, $\psi(\m)=\m$ which implies that $\pi_{\m,\psi(\m)}=\id$, so that $\pi_{\m,\psi(\m)}(\psi(\Gamma))\subseteq \Gamma$ becomes $\psi(\Gamma)\subseteq\Gamma$. Since $\#\psi(\Gamma)=\#\Gamma$, we must have $\psi(\Gamma)=\Gamma$. Thus, we have $R_{\psi(\m),\psi(\Gamma)}=R_{\m,\Gamma}$. 
\setlength{\parindent}{0cm} \setlength{\parskip}{0.5cm}

(ii): It is enough to show that $R_{\m,\Gamma}$ contains a non-unit $a$ because then $\bigcap_{n\geq 0} a^nR=\{0\}$. Observe that $R_{\m,\Gamma}$ contains the set $(1+\m_0)_+$ of totally positive elements in $1+\m_0$. Since $(1+\m_0)_+$ contains infinitely many positive integers, we see that $R_{\m,\Gamma}$ contains non-units.
\eproof
\setlength{\parindent}{0cm} \setlength{\parskip}{0.5cm}

In this setting, we have the following complete rigidity theorem:

\btheo
\label{thm:congr}
For $i=1,2$, let $K_i$ be a number field with ring of integers $R_i$, and suppose $(R_i)_{\m_i,\Gamma_i}\subseteq R_i^\times$ is a congruence monoid as in \cite[\S~3]{Bru1}. The following statements are equivalent:
\setlength{\parindent}{0cm} \setlength{\parskip}{0cm}

\begin{enumerate}[\upshape(i)]
	\item the algebraic actions $(R_1)_{\m_1,\Gamma_1}\acts R_1$ and $(R_2)_{\m_2,\Gamma_2}\acts R_2$ are isomorphic;
	\item $((R_1)_{\m_1,\Gamma_1}^{-1}R_1\rtimes \langle (R_1)_{\m_1,\Gamma_1}\rangle)\ltimes\ol{R_1}$ and $((R_2)_{\m_2,\Gamma_2}^{-1} R_2\rtimes\langle (R_2)_{\m_2,\Gamma_2}\rangle)\ltimes \ol{R_2}$ are isomorphic as topological groupoids;
    \item $(R_2)_{\m_2,\Gamma_2}^{-1} R_2\rtimes\langle (R_2)_{\m_2,\Gamma_2}\rangle \acts \ol{R_1}$ and $(R_2)_{\m_2,\Gamma_2}^{-1} R_2\rtimes\langle (R_2)_{\m_2,\Gamma_2}\rangle\reg \acts \ol{R_2}$ are continuously orbit equivalent in the sense of \cite{Li:ETDS,Li:IMRN2017};
    \item we have an isomorphism of Cartan pairs
$$
  (\fA_{(R_1)_{\m_1,\Gamma_1}\acts R_1},\fD_{(R_1)_{\m_1,\Gamma_1}\acts R_1}) \cong (\fA_{(R_2)_{\m_2,\Gamma_2}\acts R_2},\fD_{(R_2)_{\m_2,\Gamma_2}\acts R_2});
$$
    \item we have an isomorphism of Cartan pairs 
$$
\resizebox{\hsize}{!}{
  $(C_\lambda^*(R_1\rtimes (R_1)_{\m_1,\Gamma_1}),D_\lambda(R_1\rtimes (R_1)_{\m_1,\Gamma_1})) \cong (C_\lambda^*(R_2\rtimes (R_2)_{\m_2,\Gamma_2}),D_\lambda(R_2\rtimes (R_2)_{\m_2,\Gamma_2}));$
  }
$$
    \item we have an isomorphism of abstract groups
$$
  \bm{F}(((R_1)_{\m_1,\Gamma_1}^{-1}R_1\rtimes \langle (R_1)_{\m_1,\Gamma_1}\rangle)\ltimes\ol{R_1}) \cong \bm{F}(((R_2)_{\m_2,\Gamma_2}^{-1} R_2\rtimes\langle (R_2)_{\m_2,\Gamma_2}\rangle)\ltimes \ol{R_2});
$$
    \item we have an isomorphism of abstract groups
$$
  \bm{D}(((R_1)_{\m_1,\Gamma_1}^{-1}R_1\rtimes \langle (R_1)_{\m_1,\Gamma_1}\rangle)\ltimes\ol{R_1}) \cong \bm{D}(((R_2)_{\m_2,\Gamma_2}^{-1} R_2\rtimes\langle (R_2)_{\m_2,\Gamma_2}\rangle)\ltimes \ol{R_2}).
$$
\end{enumerate}
\etheo
\setlength{\parindent}{0cm} \setlength{\parskip}{0cm}

\bproof
(i)$\Leftrightarrow$(ii): Let $1\neq a\in R_{\m_i,\Gamma_i}$. Then multiplication by $1-a$ is injective on $K_i=\Qz\otimes_\Zz R_i$, and thus also injective on $R_{\m_i,\Gamma_i}^{-1}R_i\subseteq K_i$. 
Thus, by Lemma~\ref{lem:congrumonoid2}(i) and Lemma~\ref{lem:congrumonoid1}, we can apply Corollary~\ref{cor:Mispan(Mi),Oi} to obtain the desired equivalence.
\setlength{\parindent}{0cm} \setlength{\parskip}{0.5cm}

Equivalence of (ii), (iii), and (iv) follows from \cite[Theorem~2.7]{Li:IMRN2017}.

(v)$\Rightarrow$(iv) follows from the description of the primitive ideals of $C_\lambda^*(R_i\rtimes (R_i)_{\m_i,\Gamma_i})$ in \cite[Theorem~7.1]{Bru1} and the observation that $\fA_{(R_i)_{\m_i,\Gamma_i}\acts R_i}$ is the unique simple quotient of $C_\lambda^*(R_i\rtimes (R_i)_{\m_i,\Gamma_i})$  and the quotient map $C_\lambda^*(R_i\rtimes (R_i)_{\m_i,\Gamma_i})\to \fA_{(R_i)_{\m_i,\Gamma_i}\acts R_i}$ carries $D_\lambda(R_i\rtimes (R_i)_{\m_i,\Gamma_i})$ onto $\fD_{(R_i)_{\m_i,\Gamma_i}\acts R_i}$ (compare \cite[\S~8]{Bru1}). Clearly, (i)$\Rightarrow$(iv).

Equivalence of (vi), (vii), and (ii) follows from Lemma~\ref{lem:congrumonoid2}(ii) and Remark~\ref{rem:GPDCartanTFG}. 
\eproof
\setlength{\parindent}{0cm} \setlength{\parskip}{0.5cm}

Specalizing to the case where the moduli are trivial, and observing that the algebraic actions $R_1\reg\acts R_1$ and $R_2\reg\acts R_2$ are isomorphic if and only if $K_1\cong K_2$, we obtain:
\bcor
\label{cor:NT}
For $i=1,2$, let $K_i$ be a number field with rings of integers $R_i$, denote the corresponding ring C*-algebras by $\fA[R_i]$, and their canonical Cartan subalgebras by $\fD[R_i]$. The following are equivalent:
\setlength{\parindent}{0cm} \setlength{\parskip}{0cm}

\begin{enumerate}[\upshape(i)]
    \item $K_1$ and $K_2$ are isomorphic;    
    \item $(K_1 \rtimes K_1\reg) \ltimes \ol{R_1}$ and $(K_2 \rtimes K_2\reg) \ltimes \ol{R_2}$ are isomorphic as topological groupoids;
    \item $K_1 \rtimes K_1\reg \acts \ol{R_1}$ and $K_2 \rtimes K_2\reg \acts \ol{R_2}$ are continuously orbit equivalent in the sense of \cite{Li:ETDS,Li:IMRN2017};
    \item $(\fA[R_1],\fD[R_1])$ and $(\fA[R_2],\fD[R_2])$ are isomorphic as Cartan pairs;    
    \item $(C_\lambda^*(R_1\rtimes R_1\reg),D_\lambda(R_1\rtimes R_1\reg))$ and $(C_\lambda^*(R_2\rtimes R_2\reg),D_\lambda(R_2\rtimes R_2\reg))$ are isomorphic as Cartan pairs;
    \item $\bm{F}((K_1 \rtimes K_1\reg) \ltimes \ol{R_1})$ and $\bm{F}((K_2 \rtimes K_2\reg) \ltimes \ol{R_2})$ are isomorphic as abstract groups;
    \item $\bm{D}((K_1 \rtimes K_1\reg) \ltimes \ol{R_1})$ and $\bm{D}((K_2 \rtimes K_2\reg) \ltimes \ol{R_2})$ are isomorphic as abstract groups.
\end{enumerate}
\ecor
\setlength{\parindent}{0cm} \setlength{\parskip}{0.5cm}

\bremark
\label{rmk:resolve}
The equivalence of (i) and (v) in Theorem~\ref{thm:congr} completely answers the natural problem left open in \cite[\S~5.2]{BL}: The Cartan pair $(C_\lambda^*(R\rtimes R_{\m,\Gamma}),D_\lambda(R\rtimes R_{\m,\Gamma}))$ remembers the isomorphism class of the semigroup $R\rtimes R_{\m,\Gamma}$. 
The equivalences of (i), (iii), and (v) in Corollary~\ref{cor:NT}, completely answers the natural question left open in \cite[\S~1]{Li:Adv2016}.
\eremark

\subsubsection{Algebraic actions from commutative algebra}
\label{sss:commalg}
In this subsection, we analyze a class of algebraic $\Nz^d$-actions that are irreversible analogues of the algebraic $\Zz^d$-actions studied in \cite[Chapter~II]{Sch}.

We need the following observation. If $R$ is a commutative finitely generated torsion-free ring of rank $n$, then integral closure $\cO$ of $\Zz$ in $\Qz R$ is then a ring by \cite[Corollary~1.11]{Reiner}. Since $R\subseteq \cO$, for each element $a\in R$, the map $\dot{\sigma}_a\colon \Qz R\to\Qz R$, $\dot{\sigma}_a(x)=ax$, leaves $\cO$ invariant. Thus, $N(a)=|\det(\dot{\sigma}_a)|$ lies in $\Zz_{>0}$ for every $a\in R^\times$.

\blemma
\label{lem:Zkaction}
Let $R$ be a commutative finitely generated torsion-free ring of rank $n$ and $a_1,...,a_k\in R^\times\setminus R^*$. In addition, assume that for every $1\leq i\leq k$, there exists a prime $p$ such that $p\mid N(a_i)$ and $p\nmid N(a_j)$ for $j\neq i$. 
\setlength{\parindent}{0cm} \setlength{\parskip}{0cm}

\begin{enumerate}[\upshape(i)]
	\item If $\psi\in\Aut_{\Qz\textup{-alg}}(\Qz R)$ is such that $\psi(\langle a_1,...,a_k\rangle^+)\subseteq \langle a_1,...,a_k\rangle^+$, then $\psi=\id$. 
	\item $\langle a_1,...,a_k\rangle\cap\cO=\langle a_1,...,a_k\rangle^+$.
	\item $a_1,...,a_k$ are multiplicatively independent, so that the canonical map $\Nz^k\twoheadrightarrow \langle a_1,...,a_k\rangle^+$ is an isomorphism.
\end{enumerate}
\elemma
\setlength{\parindent}{0cm} \setlength{\parskip}{0cm}

\bproof
(i): Take $\psi\in\Aut_{\Qz\textup{-alg}}(\Qz R)$ with $\psi(\langle a_1,...,a_k\rangle^+)\subseteq \langle a_1,...,a_k\rangle^+$. Then for each $1\leq i\leq k$, there exist $n_1,...,n_k\in\Nz$ such that $\psi(a_i)=a_1^{n_1}\cdots a_k^{n_k}$. Now, for $a\in R^\times$, we have $\psi\circ \sigma_a=\sigma_{\psi(a)}\circ\psi$. In particular, $\det(\sigma_{\psi(a)})=\det(\sigma_a)$. Thus $N(a_i)=N(\psi(a_i))=N(a_1)^{n_1}\cdots N(a_k)^{n_k}$. This, together with  our assumption on $N(a_i)$, shows that $\psi(a_i)=a_i$.
\setlength{\parindent}{0cm} \setlength{\parskip}{0.5cm}

(ii): 
Let $x\in \langle a_1,...,a_k\rangle\cap\cO$. Then there exists $n_1,...,n_k\in\Zz$ with $x=a_1^{n_1}\cdots a_k^{n_k}$. We need to show that each $n_i$ is non-negative.
By assumption, for each $1\leq i\leq k$, there exists a rational prime $p$ dividing $N(a_i)$ such that $p\nmid N(a_j)$ for $j\neq i$. Thus, $0\leq v_p(N(x))=\sum_{j=1}^kn_jv_p(N(a_j))=n_iv_p(N(a_i))$, which shows $n_i\geq 0$ (here, the first inequality uses that $x$ lies in $\cO$).
The containment ``$\supseteq$'' is obvious. 

(iii): Suppose $a_1^{n_1}\cdots a_k^{n_k}=1$ for some $n_1,...,n_k\in\Zz$. Fix $1\leq i\leq k$, and choose a rational prime $p$ such that $v_p(N(a_i))>0$ and $v_p(N(a_j))=0$ for $j\neq i$. Now $0=v_p(N(a_1)^{n_1}\cdots N(a_k)^{n_k})=n_iv_p(N(a_i))$, so that $n_i=0$.
\eproof
\setlength{\parindent}{0cm} \setlength{\parskip}{0.5cm}

Let $d\in\Zz_{>0}$ and denote by $R_d^+:=\Zz[u_1,...,u_d]$ the ring of polynomials with integer coefficients in the $d$ variables $u_1,...,u_d$. Let $I\unlhd R_d^+$ be a non-zero ideal. By the Hilbert Basis Theorem (see, for instance, \cite[Theorem~1.2]{Eis}), $R_d^+$ is Noetherian, so that  there exists $f_1,...,f_m\in R_d^+$ such that $I$ is generated by $\{f_1,...,f_m\}$. Since we are only interested in the quotient ring $R_d^+/I$, let us assume that $u_i\notin I$ for all $1\leq i\leq d$.
Let 
\[
V(I):=\{\bm{z}\in\Cz^d : f(\bm{z})=0\text{ for every } f\in I\}\subseteq \Cz^d
\]
be the complex variety defined by $I$. It follows from \cite[Chapter~5,~Theorem~6]{CLO2} that $\dim_\Qz \Qz\otimes_{\Zz} R_d^+/I<\infty$ if and only if $V(I)$ is a finite set. If $\#V(I)<\infty$, then $\Cz\otimes_\Zz I$ is said to be zero-dimensional, in which case  there exists a basis for $\Cz\otimes_{\Zz} R_d^+/I$ consisting of (cosets of) monomials (see \cite[Chapter~5,~Proposition~4]{CLO2}), so that $R_d^+/I$ is finitely generated. For the remainder of this section, we shall assume $\#V(I)<\infty$.

For $f\in R_d^+$, let $\sigma_f$ denote the endomorphism of $R$ given by left multiplication with the coset $f+I$. Let $\chi_f(t)$ denote the characteristic polynomial of $\sigma_f$ viewed as an endomorphism of $\Cz\otimes_\Zz R_d^+$.  Let us record the following properties of these endomorphisms:
\blemma
\label{lem:mf}
For $f\in R_d^+$, we have
\setlength{\parindent}{0cm} \setlength{\parskip}{0cm}

\begin{enumerate}[\upshape(i)]
	\item $\chi_f(t)=\prod_{\bm{z}\in V(I)}(t-f(\bm{z}))^{\mu(\bm{z})}$, where $\mu(\bm{z}):=\dim_\Cz \cO_{\bm{z}}/(\Cz\otimes_\Zz I) \cO_{\bm{z}}$, where $\cO_{\bm{z}}$ is the localisation of $\Cz\otimes_\Zz R$ at the maximal ideal $\m_{\bm{z}}:=\{g \in \Cz\otimes_\Zz R : g(\bm{z})=0\}$;
	\item $\sigma_f$ is injective if and only if $f(\bm{z})\neq 0$ for every $\bm{z}\in V(I)$;
	\item $\id-\sigma_f=\sigma_{1-f}$ is injective if and only if $f(\bm{z})\neq 1$ for every $\bm{z}\in V(I)$;
	\item if for all $F\subseteq V(I)$, $\prod_{\bm{z}\in F}f(\bm{z})\neq \pm 1$, then $\chi_f$ is not divisible by any unimodular polynomial.
\end{enumerate}
\elemma
\setlength{\parindent}{0cm} \setlength{\parskip}{0cm}

\begin{proof}
(i) follows from \cite[Chapter~4,~Proposition~2.7]{CLO1}, and the other parts are consequence of this.
\end{proof}

If $N\in\Zz_{>0}$, then it follows from part (ii) of Lemma~\ref{lem:mf} that the endomorphism $\sigma_N$ is injective, so we see that $R_d^+/I$ is torsion-free. 
\setlength{\parindent}{0cm} \setlength{\parskip}{0.5cm}

We let $\dot{u}_i$ denote the image of $u_i$ modulo $I$. Note that $\dot{u}_i\neq 0$ by our assumption that $u_i\notin I$. If $\bm{z}_i\neq 0$ for all $\bm{z}\in V(I)$, then $\sigma_{u_i}$ is an injective endomorphism of $R_d^+/I$ by Lemma~\ref{lem:mf}, and we obtain an algebraic action $\mon{\dot{u}_1,...,\dot{u}_d}\acts R_d^+/I$, which satisfies \eqref{FI}. Let $R_d:=\Zz[u_1^{\pm 1},...,u_d^{\pm 1}]$ be the ring of Laurent polynomials in the variables $u_1,...,u_d$; then $R_d^+/I$ embeds in $R_d/ IR_d$, and the multiplicative group $\gp{\dot{u}_1,...,\dot{u}_d}$ acts on $R_d/IR_d$ by multiplication. It is easy to see that $\gp{\dot{u}_1,...,\dot{u}_d}\acts R_d/ IR_d$ is a globalization of $\mon{\dot{u}_1,...,\dot{u}_d}\acts R_d^+/I$.
Let $\Omega_I$ denote the completion of $R_d^+/I$ with respect to the family of $\mon{\dot{u}_1,...,\dot{u}_d}$-constructible subgroups.

\btheo
\label{thm:commalg}
For $i=1,2$, let $d_i\in\Zz_{>0}$ and let $I_i$ be a non-zero ideal of $\Zz[u_1,...,u_{d_i}]$. Assume that for $i=1,2$,
\setlength{\parindent}{0cm} \setlength{\parskip}{0cm}

\begin{enumerate}[\upshape(a)]
	\item $\#V(I_i)<\infty$ and $u_k\notin I_i$ for all $k=1,...,d_i$;
	\item $\bm{z}_k\neq 0$ for every $\bm{z}\in V(I_i)$ and $k=1,...,d_i$;
	\item there exists a monomial $f$ in $u_1, \dotsc, u_{d_i}$ such that $f(\bm{z}) \neq 1$ for all $\bm{z}\in V(I_i)$;
	\item for each $1\leq j\leq d_i$, there exists a rational prime $p$ with $p\mid N(\dot{u}_j)$ and $p\nmid N(\dot{u}_k)$ for $k\neq j$.
\end{enumerate}
\setlength{\parindent}{0cm} \setlength{\parskip}{0.5cm}

Then the following statements are equivalent:
\setlength{\parindent}{0cm} \setlength{\parskip}{0cm}

\begin{enumerate}[\upshape(i)]
	\item the algebraic actions $\Nz^{d_1}\acts R_{d_1}^+/I_1$ and $\Nz^{d_2}\acts R_{d_2}^+/I_2$ are isomorphic;
	\item $(R_{d_1}/I_1R_{d_1}\rtimes \Zz^{d_1})\ltimes \Omega_{I_1}$ and $(R_{d_2}/I_2R_{d_2}\rtimes \Zz^{d_2})\ltimes \Omega_{I_2}$ are isomorphic as topological groupoids.
\end{enumerate}
\etheo
\setlength{\parindent}{0cm} \setlength{\parskip}{0cm}

\bproof
First, note that (d) implies $|N(\dot{u}_j)|>1$, so that $\dot{u}_j$ is a non-unit for all $1\leq j\leq d_i$. Conditions (a) and (b) ensure that $R_{d_i}^+/I_i$ is a finitely generated torsion-free ring and that the action $\Nz^{d_i}\acts R_{d_i}^+/I_i$ is by injective group endomorphisms whose images all have finite index. We are in the situation of Corollary~\ref{cor:Mispan(Mi),Oi}, so we only need to show that (N) and (S) are satisfied. That (N) and (S) are satisfied follows from parts (ii) and (i) of Lemma~\ref{lem:Zkaction}, respectively.
\eproof

\begin{remark}
	In the situation of Theorem~\ref{thm:commalg}, Lemma~\ref{lem:Zkaction}(iii) implies that $\Nz^{d_i}\cong \mon{\dot{u}_1,...,\dot{u}_{d_i}}$. 
\end{remark}

\begin{remark}
	Every finitely generated commutative ring of finite rank is isomorphic to a ring of the form $\Zz[u_1,...,u_d]/I$, where $I$ is zero-dimensional. However, the isomorphism will typically not be canonical, e.g., if $R$ is the ring of algebraic integers in a number field $K$, then any choice of $\Zz$-basis $\{x_1,...,x_d\}$ for $R$ gives rise to a surjective homomorphism $\Zz[u_1,...,u_d]\to R$, whose kernel must be a zero-dimensional ideal.
\end{remark}
\setlength{\parindent}{0cm} \setlength{\parskip}{0.5cm}

Let us explain two concrete example classes that are covered by Theorem~\ref{thm:commalg}.

\bex[Principal algebraic $\Nz$-actions]
A proper ideal $I\unlhd \Zz[u]$ satisfies $\#V(I)<\infty$ if and only if $I=\Zz[u]f$ for a non-constant monic polynomial $f\in\Zz[u]$ (see, e.g., \cite[Proposition~4.1.a]{Eis}). The action $\sigma_u\colon\Nz\acts \Zz[u]/\Zz[u] f$ is called a \emph{principal algebraic $\Nz$-action}. When $f$ is non-constant and monic, the cosets of $1,u,...,u^{n-1}$ form a $\Zz$-basis for $\Zz[u]/\Zz[u] f$, where $n=\deg(f)$. The matrix for $\sigma_u$ with respect to this basis is equal to the companion matrix $C_f$ of $f$, so $\sigma_u$ is injective if and only if $f(0)=\pm\det(C_f)$ is non-zero. All in all, since $V(\Zz[u] f)$ is the set of zeros of $f$, we see that $\sigma_u\colon\Nz\acts \Zz[u]/\Zz[u] f$ satisfies conditions (a)-(d) in Theorem~\ref{thm:commalg} if and if $f$ is non-contant, monic, and $|f(0)|>1$, and $f(1)\neq 0$.
\eex
\bremark
It follows from \cite[Theorem]{Krz} that $\sigma_u\colon\Nz\acts \Zz[u]/\Zz[u] f$ is exact if and only if no unimodular polynomial divides $f$ (in $\Qz[u]$). 
\eremark

\bex[Algebraic $\Nz^d$-actions defined by a point]
A special class of $\Nz^d$-actions arises from $d$-tuples of algebraic integers. These are the irreversible analogues of the algebraic $\Zz^d$-actions from \cite[\S~7]{Sch}.  Suppose that $\bm{c}=(c_1,...,c_d)\in (\overline{\Zz}^\times)^d$, where $\overline{\Zz}$ denotes the ring of all algebraic integers, and let $\p_{{\bm{c}}}$ denote the kernel of the evaluation at $\bm{c}$ map $R_d^+\to \Zz[c_1,...,c_d]\subseteq\Qz(c_1,...,c_d)$. 
Then $\p_{\bm{c}}$ is a prime ideal of $R_d^+$, and we can characterize when the action $\Nz^d\acts R_d^+/\p_{\bm{c}}$ satisfies conditions (a)-(d) in Theorem~\ref{thm:commalg} in terms of $\bm{c}$. First, identify $V(\p_{\bm{c}})$ with the set $\Hom(\Qz(c_1,...,c_d),\overline{\Qz})$ of field embeddings of $\Qz(c_1,...,c_d)$ into $\overline{\Qz}$, the algebraic closure of $\Qz$ in $\Cz$. Explicitly, this identification is given by sending $\bm{z}=(z_1,...,z_d)\in V(\p_{\bm{c}})$ to the embedding $\Qz(c_1,...,c_d)\hookrightarrow\overline{\Qz}$ determined by $c_i\mapsto z_i$.
From this, it is easy to see that conditions (a) and (b) from Theorem~\ref{thm:commalg} are satisfied if and only if $c_i\neq 0$ for all $i$, and condition (c) from Theorem~\ref{thm:commalg} is satisfied if and only if there exists a finite non-empty set $F\subseteq \{1,...,d\}$ such that $\prod_{i\in F}c_i\neq 1$; to see this, note that for any $\bm{z}=(z_1,...,z_d)\in V(\p_{\bm{c}})$, we have $\prod_{i\in F}z_i\neq 1$ if and only if $\prod_{i\in F}c_i\neq 1$.
Condition (d) from Theorem~\ref{thm:commalg} is satisfied if and only if for each $1\leq j\leq d_i$, there exists a rational prime $p$ with $p\mid N(c_j)$ and $p\nmid N(c_k)$ for $k\neq j$.
\eex

\subsection{Algebraic actions from rings: The non-commutative case}
 \label{s:ncrings}

For $i=1,2$, let $R_i$ be a ring whose additive group is finitely generated and torsion-free, let $\cL_i\subseteq R_i$ be full rank subgroup, and let $M_i\subseteq R_i\reg$ a submonoid such that $\cL_i$ is $M_i$-invariant. Then, $M_i\acts \cL_i$ is faithful since $\cL_i$ has finite index in $R_i$ and $R_i$ is torsion-free. Let $\ol{\cL}_i$ denote the completion of $\cL_i$ with respect to the family $\cC_i$ of $M_i$-constructible subgroups of $\cL_i$. Our goal now is to establish the following rigidity result:
\btheo
\label{thm:ssQ}
Continue with the notation and assumptions above, with the additional assumptions that, for $i=1,2$, $\spn_\Zz(M_i)$ has finite index in $R_i$, that there exists $\kappa_i \in M_i$ for some $\kappa_i \in \Zz \setminus \gekl{0,1}$, and that $\Qz R_i$ is a semisimple $\Qz$-algebra, i.e., the (Jacobson) radical of $\Qz R_i$ is trivial (see, e.g., \cite[Part II,~\S~1]{Kap}). If there is an isomorphism of topological groupoids
\[
(\Qz R_1\rtimes \langle M_1\rangle)\ltimes\ol{\cL}_1\cong (\Qz R_2\rtimes \langle M_2\rangle)\ltimes\ol{\cL}_2,
\]
then there are $\Qz$-algebra isomorphisms $\varphi_1\colon\Qz R_1\overset{\cong}{\to} \Qz R_2$ and $\varphi_2\colon\Qz R_2\overset{\cong}{\to} \Qz R_1$ such that $\varphi_1(\langle M_1\rangle)\subseteq \langle M_2\rangle$ and $\varphi_2(\langle M_2\rangle)\subseteq \langle M_1\rangle$.
\etheo
\setlength{\parindent}{0cm} \setlength{\parskip}{0cm}

Taking $M_i= R_i^\times$ and $\cL_i= R_i$ in Theorem~\ref{thm:ssQ} yields the following:
\bcor
\label{cor:ssQ}
If $\Qz R_1$ and $\Qz R_2$ are semisimple $\Qz$-algebras and there is an isomorphism of topological groupoids 
\[
(\Qz R_1\rtimes (\Qz R_1)^*)\ltimes\ol{R}_1\cong (\Qz R_2\rtimes (\Qz R_2)^*)\ltimes\ol{R}_2,
\]
then there is a $\Qz$-algebra isomorphism $\Qz R_1\cong \Qz R_2$.
\ecor
\bproof
Observe that $\spn_\Zz(R_i^\times) =  R_i$ as shown in Theorem~\ref{thm:Bhargava}. The claim now follows from Theorem~\ref{thm:ssQ}.
\eproof

\bremark
The actions $R_i^\times\acts  R_i$ in Corollary~\ref{cor:ssQ} are exact, so Remark~\ref{rem:GPDCartanTFG} applies here.
\eremark
\setlength{\parindent}{0cm} \setlength{\parskip}{0.5cm}

Before proceeding to the proof, let us explain several example classes to which our results apply.

\bex[Matrices over orders in number fields]
\label{ex:MnK}
Let $n\in\Zz_{>0}$. Let $R$ be an order in a number field $K$, and let $I\unlhd R$ be a nonzero ideal. Then, $\M_n(I)\subseteq \M_n(R)$
is invariant under the canonical action $\M_n(R)^\times\acts \M_n(R)$, so we get an algebraic action $\M_n(R)^\times\acts \M_n(I)$. We have $\Qz \M_n(I)=\M_n(K)$, so $\M_n(I)$ has full rank in $\M_n(R)$. We have $\Zz 1_n\subseteq \M_n(R)^\times$, and $\spn_\Zz(\M_n(R)^\times)$ has finite index in $\M_n(R)$ by the proof of Corollary~\ref{cor:ssQ}. Thus, the hypotheses of Theorem~\ref{thm:ssQ} are satisfied (for $\cL = \M_n(I)$ and $M = \M_n(R)\reg$). Moreover, $\gp{\M_n(R)^\times}=\GL_n(K)$. Therefore, if $K_1$ and $K_2$ are number fields with rings of algebraic integers $R_1$ and $R_2$, respectively, $I_1\unlhd R_1$, $I_2\unlhd R_2$ are non-zero ideals, $n_1,n_2\in\Zz_{>0}$, and if $(\M_{n_1}(K_1) \rtimes \GL_{n_1}(K_1))\ltimes\ol{\M_{n_1}(I_1)} \cong (\M_{n_2}(K_2) \rtimes \GL_{n_2}(K_2))\ltimes\ol{\M_{n_2}(I_2)}$ as topological groupoids, then $\M_{n_1}(K_1) \cong \M_{n_2}(K_2)$.
\setlength{\parindent}{0.5cm} \setlength{\parskip}{0cm}

In particular, Theorem~\ref{thm:ssQ} implies that 
\[
(\M_{n_1}(K)\rtimes\GL_{n_1}(K))\ltimes \M_{n_1}(\ol{R_1}) \cong (\M_{n_2}(K_2)\rtimes\GL_{n_2}(K_2))\ltimes \M_{n_2}(\ol{R_2})
\]
as topological groupoids if and only if $n_1=n_2$ and $K_1\cong K_2$.
\eex
\setlength{\parindent}{0cm} \setlength{\parskip}{0.5cm}

\bex[Group rings of finite groups]
\label{ex:grprings}
Let $F_1$ and $F_2$ be finite groups. By Maschke's Theorem (see, e.g., \cite[Theorem~25]{Kap}), $\Qz F_i$ is a semisimple $\Qz$-algebra ($i=1,2$), so Corollary~\ref{cor:ssQ} implies the following: If there is an isomorphism $(\Qz F_1\rtimes (\Qz F_1)^*)\ltimes\ol{\Zz F_1}\cong (\Qz F_2\rtimes (\Qz F_2)^*)\ltimes\ol{\Zz F_2}$, then $\Qz F_1\cong \Qz F_2$. Note that if $F_1,F_2$ are Abelian, then $\Qz F_1\cong \Qz F_2$ if and only if $F_1\cong F_2$ by \cite[Corollary~1~\&~Theorem~3]{PW}. It is a non-trivial result that there exist finite non-Abelian groups $F_1,F_2$ with $\Zz F_1\cong\Zz F_2$ and $F_1\not\cong F_2$ (see \cite{Hertweck}).
\eex

\bex[Central simple algebras over number fields]
\label{ex:ordersinssalgs}
Let $\cA$ be a central simple algebra over the number field $K$, i.e., $\cA$ is a finite-dimensional simple $K$-algebra whose centre is precisely $K$, and let $\cO$ be an order in $\cA$. By the Wedderburn Structure Theorem, there exists a (central) division algebra $D$ over $K$ and $\mu \in\Zz_{>0}$ such that $\cA\cong \M_{\mu}(D)$. Thus, the algebraic action $\cO^\times\acts \cO$ fits into the setting of Corollary~\ref{cor:ssQ}. Thus if $(\cA_1 \rtimes \cA_1^*) \ltimes \ol{\cO}_1 \cong (\cA_2 \rtimes \cA_2^*) \ltimes \ol{\cO}_2$ as topological groupoids, then $\cA_1 \cong \cA_2$.
\eex

Now our goal is to prove Theorem~\ref{thm:ssQ}. For the remainder of this section, we shall work with the assumptions and notation from Theorem~\ref{thm:ssQ}. We need some preparations.

\blemma
\label{lem:nil}
Let $R$ be a ring whose additive group is torsion-free and of rank $n$, and let $M\subseteq  R^\times$ be a submonoid. Suppose that $\alpha\in \langle M\rangle$ satisfies $\alpha=\gamma\alpha^{\kappa} \gamma^{-1}$ for some $\gamma\in (\Qz R)^*$ and $\kappa \in\Zz_{>1}$. Set $m:= \kappa (\dim_\Qz\Qz  R)!-1$. Then there exists a nilpotent element $\eta_\alpha\in \Qz R$ such that $\alpha^m=1+\eta_\alpha$ (i.e., $\alpha^m$ is a unipotent element of the $\Qz$-algebra $\Qz R$).
\elemma
\setlength{\parindent}{0cm} \setlength{\parskip}{0cm}

\bproof
The map $\pi\colon \Qz R \to \End_\Cz(\Cz R)\cong \M_n(\Cz)$ given by $\pi(q\otimes a)(z\otimes b)=qz\otimes ab$ is an injective $\Qz$-algebra homomorphism. The equation $\alpha=\gamma\alpha^{\kappa} \gamma^{-1}$ implies that $\pi(\alpha)=\pi(\gamma)\pi(\alpha)^{\kappa} \pi(\gamma)^{-1}$ in $\M_n(\Cz)$. It follows that $\sp(\pi(\alpha))=\sp(\pi(\alpha)^{\kappa})=\sp(\pi(\alpha))^{\kappa} :=\{\lambda^{\kappa} : \lambda\in\sp(\pi(\alpha))\}$, where $\sp(\pi(\alpha))$ is the spectrum of $\pi(\alpha)$. Thus, the map $\sp(\pi(\alpha))\to \sp(\pi(\alpha))$ given by $\lambda\mapsto \lambda^k$ is bijective; write $\sp(\pi(\alpha))=\{\lambda_1,...,\lambda_j\}$, and let $\rho$ be the permutation of $\sp(\pi(\alpha))$ determined by $\lambda_i=\lambda_{\rho(i)}^{\kappa}$ for all $1\leq i\leq j$. Since $j\leq \dim_\Qz\Qz R$, we have $\rho^{(\dim_\Qz\Qz R)!}=\id$. Thus, $\lambda_i=\lambda_{\rho^{(\dim_\Qz\Qz R)!}(i)}^{\kappa (\dim_\Qz\Qz R!)}=\lambda_i^{\kappa (\dim_\Qz\Qz R)!}$, so that $\lambda_i^{\kappa (\dim_\Qz\Qz R)!-1}=1$. We now see that $1$ is the only eigenvalue of $\pi(\alpha)^m$. By considering the Jordan Normal Form of $\pi(\alpha)^m$, it follows that there exists a nilpotent matrix $N_\alpha\in \M_n(\Cz)$ such that $\pi(\alpha)^m=1+N_\alpha$. Since $N_\alpha=\pi(\alpha)^m-1=\pi(\alpha^m-1)$, we see, by injectivity of $\pi$, that $\eta_\alpha:=\alpha^m-1$ is nilpotent. 
\eproof
\setlength{\parindent}{0cm} \setlength{\parskip}{0.5cm}

Given a division algebra $D$, we shall regard $D^n$ as an $\M_n(D)$-$D$-bimodule in the usual way. Thus, a basis for $D^n$ will always mean a right $D$-basis. We shall use a subscript $D$ on the right to indicate that we are viewing something as a right $D$-vector space.

\blemma
\label{lem:rkGamma}
Let $D$ be a finite dimensional division algebra over $\Qz$ and $n\in\Zz_{>0}$. Suppose that $\Sigma$ is a non-trivial finitely generated Abelian subgroup of $\GL_n(D)$ consisting of unipotent matrices, so that every $\alpha\in\Sigma$ is of the form $\alpha=1+\eta_\alpha$, where $\eta_\alpha\in \M_n(D)$ is a nilpotent matrix. Let $\mfk:=\bigcap_{\alpha\in\Sigma}\ker\eta_\alpha$ and $k:=\dim\mfk_D$. We have $\rk_\Zz\Sigma< n\cdot(n-k)[D:\Qz]$.
\elemma
\setlength{\parindent}{0cm} \setlength{\parskip}{0cm}

\bproof
Let $\fN:=\spn_\Qz\{\eta_\alpha : \alpha\in\Sigma\} \subseteq \M_n(D)$. For $\alpha,\alpha'\in\Sigma$, we have 
\[
1+\eta_{\alpha\alpha'}=\alpha\alpha'=(1+\eta_\alpha)(1+\eta_{\alpha'})=1+\eta_\alpha+\eta_{\alpha'}+\eta_{\alpha}\eta_{\alpha'},
\]
so that $\eta_{\alpha}\eta_{\alpha'}=\eta_{\alpha\alpha'}-\eta_\alpha-\eta_{\alpha'}$ lies in $\fN$. From this, we see that $\fN$ is a non-unital commutative sub-$\Qz$-algebra of $\M_n(D)$. 
\setlength{\parindent}{0cm} \setlength{\parskip}{0.5cm}

For $\alpha\in\Sigma$, let $\log\alpha:=\sum_{i=1}^\infty(-1)^{i-1}\frac{(1-\alpha)^i}{i}$; this is a finite sum because $1-\alpha$ is nilpotent. Note that $\log\alpha$ lies in $\fN$. Since $\Sigma$ is Abelian, $\log$ defines an injective group homomorphism (with inverse given by $\exp$, see for instance \cite[\S~2.10, Exercise~8]{Hall}) from $\Sigma$ into $\fN$, so that $\rk_\Zz\Sigma\leq \dim_\Qz\fN$. Thus, we will be done once we show that $\dim_\Qz\fN<n\cdot(n-k)[D:\Qz]$.

Since $\fN$ is closed under multiplication and consists of nilpotent elements, \cite[Part~II,~\S~5,~Theorem~35]{Kap} asserts that there exists a right $D$-basis $w_1,...,w_n$ for $D^n$ such that $\fN$ is strictly upper triangular with respect to this basis; this means that if we define $W_i:=\spn\{w_1,...,w_i\}_D$ for $i=1,..,n$, then $\eta w_1=0$ and for each $i\geq 2$, $\eta W_i\subseteq W_{i-1}$ for all $\eta\in\fN$. Let $\cW\subseteq \{w_1,...,w_n\}$ be a subset such that $\cW$ together with a basis for $\mfk$ is a basis for $D^n$, and put $W = \lspan(\cW)_D$. 
Since $\eta\vert_\mfk\equiv 0$ for all $\eta\in\fN$, the map $\fN \to \Hom(W,D^n)_D, \, \eta\mapsto \eta\vert_W$ is injective; moreover, since $w_n\not\in\fN D^n$, we see that it is not surjective. Hence, $\dim_\Qz\fN<\dim_{\Qz}\Hom(W,D^n)_D=\dim_\Qz M_{n\times(n-k)}(D)=n\cdot (n-k)[D:\Qz]$.
\eproof
\setlength{\parindent}{0cm} \setlength{\parskip}{0.5cm}

\bproof[Proof of Theorem~\ref{thm:ssQ}]

Let $\mfc$ be the cocycle defined by the composition
\[
(\Qz R_1\rtimes\langle M_1\rangle)\ltimes\ol{\cL}_1\cong (\Qz R_2\rtimes \langle M_2\rangle)\ltimes\ol{\cL}_2\overset{(h,\bm{y})\mapsto h}{\longrightarrow} \Qz R_2\rtimes\langle M_2\rangle,
\]
where the second map is the canonical cocycle obtained by projecting onto the group component. Lemma~\ref{lem:fisubgroup} produces a (finite index) subgroup $C \in\cC_1$ such that $\mfg(x):= \mfc(x,\bm{0})$ defines an injective group homomorphism from $C$ into $\Qz R_2\rtimes\langle M_2\rangle$. Let $\fT \subseteq \langle M_2\rangle$ be the image of $C$ under the composition 
\begin{equation}
\label{e:C-QRM-M}
    C \overset{\mfg}{\to} \Qz R_2\rtimes\langle M_2\rangle\overset{(b,t)\mapsto t}\to\langle M_2\rangle.
\end{equation}
By assumption, there exists $\kappa \in M_1$ with $\kappa \in \Zz \setminus \gekl{0,1}$. Let $m:= \kappa (\dim_\Qz\Qz R_2)!-1$. Then we claim that for every $\alpha\in\fT^m$, there exists a nilpotent element $\eta_\alpha\in\Qz R_2$ such that $\alpha=1+\eta_\alpha$. Indeed, as $\kappa \in M_2$, Proposition~\ref{prop:GammaAlphaEq} implies that there exists $\gamma\in\langle M_2\rangle$ such that $\alpha=\gamma\alpha^{\kappa} \gamma^{-1}$ for all $\alpha\in\fT$. The result now follows from Lemma~\ref{lem:nil}. In particular, $\fT^m$ is torsion-free. The subgroup $m C \subseteq C$ is mapped onto $\fT^m$ under the projection in \eqref{e:C-QRM-M}, so $\fT^m$ is moreover finitely generated. We will show that $\fT^m$ is trivial, which will imply that there exists an injective group homomorphism $\mfb\colon mC\to \Qz R_2$ such that $\mfg(x)=(\mfb(x),1)$ for all $x\in mC$. Since $\fT^m$ is free abelian, there exists a subgroup $C_{\fT}\subseteq mC$ such that $mC=C' \oplus C_{\fT}$, where $C'=\{x \in mC : \mfg(x) = (y,1) \text{ for some } y \in \Qz R_2\}$, and $C_{\fT}$ is mapped isomorphically onto $\fT^m$ under the map in \eqref{e:C-QRM-M}.

Let $\fB$ be the be the image of $C'$ under the composition $mC \overset{\mfg}{\to} \Qz R_2\rtimes\langle M_2\rangle\overset{(b,t) \mapsto b}{\to} \Qz R_2$. Note that $\fB$ is a subgroup of $\Qz R_2$ because the second map is a homomorphism on $\Qz R_2\rtimes\{1\}$, which contains the image of $C'$.  

Let us now show that the group $\fT^m$ is trivial. Since $\Qz R_2$ is a semisimple $\Qz$-algebra, the Artin--Wedderburn theorem implies that there exists a decomposition of $\Qz$-algebras $\Qz R_2=\bigoplus_{i=1}^r\M_{n_i}(D_i)$, where $r,n_1,...,n_r\in\Zz_{>0}$ and each $D_i$ is a finite-dimensional division algebra over $\Qz$. Thus, $\langle M_2\rangle\subseteq \bigoplus_{i=1}^r\GL_{n_i}(D_i)$. For each $i=1,...,r$, let $\fB_i$ be the image of $\fB$ under the canonical projection $\Qz  R_2\to \M_{n_i}(D_i)$. For each $i=1,...,r$, let $\fT_i$ be the image of $\fT^m$ under the canonical projection $\langle M_2\rangle\to \GL_{n_i}(D_i)$. Then $\fT_i=\{1+\eta_\alpha^i : \alpha\in\fT^m\}$, where $\eta_\alpha^i$ denotes the $i$-th coordinate of $\eta_\alpha$ (which come from Lemma~\ref{lem:nil}); in particular, the group $\fT_i$ consists of unipotent elements. Let $\k_i:=\bigcap_{\alpha\in\fT^m}\ker(\eta_\alpha^i)$ and $k_i \defeq \dim (\k_i)_D$.

Take $x'\in C'$ and write $\mfg(x')=(\beta',1)$. If $x\in C_\fT$, we can write $\mfg(x)=(\beta,\alpha)$. Then $\mfg(x') \mfg(x)=(\beta'+\beta,\alpha)$, whereas $\mfg(x) \mfg(x')=(\beta+\alpha\beta',\alpha)$. Thus, $\alpha\beta'=\beta'$, i.e., $\eta_\alpha\beta'=0$. 
It follows that $\eta_\alpha^i\beta=0$ for all $\alpha\in\fT^m$ and $\beta\in\fB_i$, so that $\im(\beta)\subseteq\k_i$ for all $\beta\in \fB_i$. From this, we see that $\im(\beta)\subseteq\k_i$ for all $\beta\in\spn(\fB_i)_{D_i}:=\{\sum_j \beta_jd_j : \beta_j\in\fB_i, d_j\in D_i\}$. Hence, $\dim(\spn(\fB_i)_{D_i})_{D_i}\leq n_i\cdot k_i$.
Now we have $\rk_\Zz\fB_i=\dim_\Qz(\Qz\otimes\fB_i)\leq \dim_\Qz(\spn(\fB_i)_{D_i})\leq n_i\cdot k_i\cdot[D_i:\Qz]$.

Assume for the sake of contradiction that $\fT^m$ is non-trivial, so that $\fT_i$ is non-trivial for some $i$ (i.e., $k_i<n_i$). By Lemma~\ref{lem:rkGamma}, we have $\rk_\Zz\fT_i<n_i\cdot (n_i-k_i)[D_i:\Qz]$, where $k_i:=\dim(\k_i)_{D_i}$. Since $\rk_\Zz\fT^m\leq\sum_{i=1}^r\rk_\Zz\fT_i$ and $\rk_\Zz\fB \leq\sum_{i=1}^r\rk_\Zz \fB_i$, we obtain
\begin{align*}
\dim_\Qz\Qz R_1 &= \rk_\Zz C=\rk_\Zz C' + \rk_\Zz C_{\fT}=\rk_\Zz \fB + \rk_\Zz \fT^m \\
&\leq \sum_{i=1}^r\rk_\Zz\fT_i+\sum_{i=1}^r\rk_\Zz\fB_i < \sum_{i=1}^r n_i^2\cdot[D_i:\Qz]=\dim_\Qz R_2,    
\end{align*}
where the strict inequality uses our assumption that $k_i<n_i$ for some $i$. By symmetry, we also get $\dim_\Qz\Qz R_2<\dim_\Qz\Qz R_1$, which is a contradiction. Thus, $k_i=n_i$ for all $i$. Hence $\fT^m$ is indeed trivial.

We conclude that there exists a group homomorphism $\mfb \colon mC \to \Qz R_2$ such that $\mfg(x)=(\mfb(x),1)$ for every $x\in mC$. Let $\mft$ be the composition 
$\langle M_1\rangle \overset{\mfg}{\to}  \Qz R_2\rtimes\langle M_2\rangle\twoheadrightarrow  \langle M_2\rangle$, where the second arrow is the canonical projection map. It is easy to check that $\mft$ is a group homomorphism, and Lemma~\ref{lem:alpha-inj} shows that $\mft$ is injective. Now Lemma~\ref{lem:equi} shows that for all $s \in M_1$ and $x \in mC$, we have $\mfb((\sigma_1)_s(x)) = (\ti{\sigma}_2)_{\mft(s)}(\mfb(x))$, where $\sigma_1$ is the algebraic action $M_1 \acts \cL_1$ and $\ti{\sigma}_2$ is the algebraic action $\gp{M_2} \acts \Qz  R_2$. Hence the same proof as for Theorem~\ref{thm:mainthm}~(i) produces an injective homomorphism $\dot{\mfb} \colon \Qz  R_1 \into \Qz  R_2$ which is equivariant. Applying Proposition~\ref{prop:extension} yields the desired result.
\eproof

\end{document}